\documentclass[a4paper]{article}
\usepackage{amsmath}
\usepackage{amsfonts,amssymb,amsthm}
\usepackage[final]{graphicx}

\newtheorem{thm}{Theorem}
\newtheorem{corr}{Corollary}
\newtheorem{lem}{Lemma}
\newtheorem{prop}{Proposition}
\newtheorem{rem}{Remark}
\newtheorem{df}{Definition}

\newcommand{\RM}{\mathbb{R}}
\newcommand{\ZM}{\mathbb{Z}}

\newcommand{\CM}{\mathbb{C}}
\newcommand{\MM}{\,\mbox{\bf M}}
\newcommand{\WM}{\,\mbox{\bf W}}

\newcommand{\HM}{\,\mbox{\bf H}}
\newcommand{\UM}{\,\mbox{\bf U}}

\newcommand{\NM}{\,\mbox{\bf N}}
\newcommand{\tr}{\,\mbox{\rm tr}}

\newcommand{\sgn}{\operatorname{sgn}}
\newcommand{\cof}{\operatorname{cof}}

\newcommand{\spec}{\operatorname{spec}}

\title{The Modulational Instability for a
Generalized Korteveg-DeVries equation}
\author{Jared C. Bronski \and Mathew A. Johnson}
\date{Submitted to ARMA on October 3, 2008}

\theoremstyle{notation}

\begin{document}
\bibliographystyle{plain}
\maketitle

\begin{abstract}
We study the spectral stability of a family of
periodic standing wave solutions to the generalized KdV (g-KdV) in a neighborhood
of the origin in the spectral plane using what amounts to a rigorous
Whitham modulation theory calculation.  In particular we are interested
in understanding the role played by the null directions of the
linearized operator in the stability of the traveling wave
to perturbations of long wavelength.

A study of the normal form of the characteristic polynomial
of the monodromy map (the periodic
Evan's function) in a neighborhood of the origin in the spectral
plane leads to two different instability indices. The first index
counts modulo $2$ the total number of periodic eigenvalues on the
real axis. This index is a generalization of the one which governs the
stability of the solitary wave. The second index provides a necessary
and sufficient condition for the existence of a long-wavelength
instability. This index is essentially the quantity calculated by
H{\v a}r{\v a}gu{\c s} and Kapitula in the small amplitude limit.
Both of these quantities can be expressed in terms of the map between
the constants of integration for the ordinary differential equation
defining the traveling waves and the conserved quantities of the
partial differential equation. These two indices together provide a
good deal of information about the spectrum of the linearized operator.

We sketch the connection of this calculation to a study of the linearized
operator - in particular we perform a perturbation calculation in terms
of the Floquet parameter.  This suggests a geometric interpretation
attached to the vanishing of the 
modulational instability index previously mentioned.

\end{abstract}

\tableofcontents

\section{Introduction and Preliminaries}

In this paper, we consider standing wave solutions to
the generalized KdV (gKdV) equation
\begin{equation}
 u_t=u_{xxx}+\left(f(u)\right)_x- c u_x\label{gkdv}
\end{equation}
where $f(\cdot)\in C^2(\RM)$ is a prescribed nonlinearity and $c$
is the wavespeed. Such solutions represent
traveling wave solutions to the generalized KdV equation
with nonlinearity $f(u)$. Of particular interest is the
case of power law nonlinearity $f(u)=u^{p+1}$, which in the cases
of $p=1,2$ represents the equations for traveling wave solutions
to the KdV and MKdV, respectively.
Obviously such traveling wave solutions are reducible to quadrature:
they satisfy
\begin{align}
u_{xx}+f(u)-c u &=a \label{dphi1} \\
\frac{u_x^2}{2} + F(u) - c\frac{u^2}{2} &= a u + E. \label{dphi2}
\end{align}
We are interested in the spectrum of the linearized operator (in the moving coordinate frame)
\[
\mu v = v_{xxx} + (v f'(u))_x - c v_x
\]
in two related settings.  First, we study the spectrum in
a neighborhood of $\lambda=0$. Physically this amounts to
long-wavelength perturbations of the underlying wave profile: in essence
slow modulations of the traveling wave. There is a well developed
physical theory, commonly known as Whitham modulation theory\cite{W2,W1}, for
dealing with such problems.
On a mathematical level the origin in the
spectral plane is distinguished by the fact that the ordinary
differential equation giving the traveling wave profile is completely
integrable. Thus the tangent space to the manifold of traveling wave profiles
can be explicitly computed, and the null-space to the linearized operator
can be built up from elements of this tangent space.
We show that these considerations give a
rigorous normal form for the spectrum of the linearized operator
in the vicinity of the origin providing that certain genericity
conditions are met.
Assuming that
these genericity conditions are met we are able to show the following:
there is a discriminant $\Delta$ which can be calculated explicitly.
If
this discriminant is positive then the spectrum in a neighborhood of the
origin consists of the imaginary axis\footnote{Note that this
does {\em not} imply spectral stability since there is the possibility
of bands of spectrum off of the imaginary axis away from the origin. } with multiplicity three.
If this discriminant is negative the spectrum of the linearization
in the neighborhood of the origin consists of the imaginary axis
(with multiplity one) together with two curves which leave the origin
along lines in the complex plane, implying instability. Long wavelength
theories are invariably geometric in nature, and the one presented here
is no exception: both the instability index and the genericity conditions
admit a natural geometric interpretation.

Secondly, we are interested in determining sufficient conditions for the
existence of unstable spectrum supported away from $\lambda=0$.  Here, this is accomplished
by calculating an orientation index using Evans function techniques: essentially comparing
the behavior of the Evans function near the origin with the asymptotic behavior
near
infinity.  Physically, such an instability amounts to an instability with respect
to finite wavelength perturbations.  The derived stability index is a generalization
of the one which governs stability of solitary waves.  In fact, in the case of power-law nonlinearity
and wave speed $c>0$, we show that in a long wavelength
limit the sign of this index, which is actually what determines stability,
 agrees with the sign of the solitary
wave stability index derived by, for example, Pego and Weinstein\cite{PW1,PW2}.

This paper uses ideas from both stability theory and modulation theory,
and thus there is an extensive background literature. Most obviously
is the stability theory of solitary wave solutions to KdV type equations
which was pioneered by Benjamin\cite{Ben} and further developed by
Bona\cite{Bona}, Grillakis\cite{G}, Grillakis,Shatah and Strauss\cite{GSS},Bona, Souganides and Strauss\cite{BSS}, Pego and Weinstein\cite{PW1,PW2}, Weinstein\cite{MIW1,MIW2} and others. In this theory the role of the discriminant is
played by the derivative of the momentum with respect to wave-speed. Our
discriminant is somewhat more complicated, which is to be expected:
the solitary waves homoclinic to the origin are a codimension two
subset of the family of periodic solutions, so one expects that the
general stability condition will more complicated. There are also a number of
calculations of the stability of periodic solutions to perturbations of
the same period, or to perturbations of twice the period, due to
Angulo Pava\cite{AP}, Angulo Pava, Bona and Scialom\cite{APBS} and others.
In this setting the linearized operator has a compact
resolvent, so the spectrum is purely discrete, and the arguments are
similar in spirit to those for the solitary wave stability. In contrast we
consider
the case of a general $L_2$ case where one must understand the continuous
spectrum of the operator.

A stability calculation in the spirit of modulation theory was given by
Rowlands\cite{Row} for the cubic nonlinear Schr\"odinger equation.
Other stability calculations in the same spirit, but differing greatly
in details and approach, were  given by Gallay and
H\v{a}r\v{a}gu\c{s}\cite{GH1,GH2}, H\v{a}r\v{a}gu\c{s} and
Kapitula\cite{HK}, Bridges and Rowlands\cite{BR}, and Bridges and Mielke\cite{BM}.
The work of Gardner\cite{Gard} is also related, though it should
be noted that the long-wavelength limit in Gardner is very different from
the one we consider here: in the former it is the traveling wave itself which
has a long period. In our calculation the period is fixed and we are
considering {\it perturbations} of long period. The current paper also
owes a debt to the substantial literature on
Whitham theory for integrable systems developed by Lax and Levermore\cite{LL1,LL2,LL3},
Flashka, Forest and McLaughlin\cite{FFM}, and many others. We note, however,
that the calculation outlined in this paper is not an integrable calculation.
The papers that are perhaps closest to that presented here are those by
Oh and Zumbrun\cite{OZI,OZII,OZIII} and Serre\cite{serre}
on the stability of periodic solutions to viscous conservation laws, where
similar results relating the behavior of the linearized spectral problem
in a neighborhood of the origin to a formal theory of slow modulations
are proved.

Our results are most explicit in the case of power law nonlinearity.
It should be noted that due to the scaling invariance in
this case  we can always
assume that $c\in\{-1,+1\}$.  Indeed, it is easy to check that if
$u(x,t;c)$ solves \eqref{gkdv} with the nonlinearity $f(u)=u^{p+1}$, then
\begin{equation}
u(x,t;c)=|c|^{1/p}u\left(|c|^{1/2}x,|c|^{3/2}t,\sgn(c)\right)\label{scaling}
\end{equation}
solves \eqref{gkdv} with wave speed $\sgn(c)$.

The paper is organized as follows: in the second section we lay out some
basic general properties of the spectrum of the linearized operator.
In the third section we explicitly compute the monodromy map
and associated periodic Evan's function at the origin. A perturbation
analysis in the neighborhood of the origin gives a normal form for the
Evan's function.
In the fourth section we develop similar results from the point of view of
the linearized operator: we compute the generalized null-space of the
linearized operator in terms of the tangent space to the ordinary
differential equation defining the traveling wave. The structure of this
null-space (under some genericity conditions) reflects that of the monodromy
map at the origin, and a similar perturbation analysis gives a normal form for the spectrum.
While the two approaches are in principle the
same most of the calculations are more easily carried out in the context of
the monodromy map/Evans function. We mainly present the calculations at
the level of the linearized operator since it helps to clarify some
aspects of the monodromy calculation.
Finally we end with some concluding remarks.

It should be noted we restrict neither the size of the periodic
solution nor the period.  Moreover, all of our analysis applies to both
localized and bounded perturbations of the underlying wave.  Also in
this paper ``stability'' will always means spectral stability.

\subsection{Preliminaries}
Note that the partial differential equation has (in general) three
conserved quantities
\begin{align}
M &= \int_0^T u(x,t) dx\nonumber \\
P &= \int_0^T u^2(x,t) dx\nonumber \\
H &= \int_0^T \frac{1}{2} u_x^2 + F(u) dx \nonumber
\end{align}
which correspond to the mass, momentum and Hamiltonian of the
solution respectively. These quantities considered as functions
of the traveling waves parameters will form an important part of the analysis.

The periodic standing wave solutions of
\eqref{gkdv} are of the form $u(x,t)= u(x)$ where $ u$ is a periodic
function in the x-variable.  Substituting this into \eqref{gkdv}
and integrating twice we see that $ u$ satisfies
\begin{equation}
\frac{1}{2} u_x^2+F\left( u\right)- \frac{c}{2}u^2-a u =E \label{dphi2}
\end{equation}
where $a,E\in\mathbb{R}$ are constants of integration and $F'=f$.
Note that the solitary wave case corresponds to $a=0,E=0,$ so the
solitary  waves are a codimension two subset of the periodic waves.
In order to assure the existence of a periodic orbits,
we must require that the effective potential
\[
V(u;a,c) = F(u) -\frac{c}{2} u^2 -au
\]
has a local minimum.
Note that this places a condition on the allowable parameter
regime $\mathcal{D}$ for our problem. We will always assume that
we are in the interior of this open region, and that the roots $u_+,u_-$
of the equation $V(u;a,c)=E$ with $V(u;a,c)<E$ for $u\in(u_{-},u_{+})$ are
simple, guaranteeing that they are $C^1$ functions of $a,E,c$.

As is standard, the the period of the corresponding
periodic orbit is given by
\begin{equation}
T=T(a,E,c):= 2\int_{ u_{-}}^{u_{+}} \frac{d
u}{\sqrt{2\left(E-V(u;a,c)\right)}}.\label{period}
\end{equation}
The above interval can be regularized at the square root branch
points $ u_-,\; u_{+}$ by the following procedure:
Write $E-V(u;a,c)=(r- u_{-})(u_{+}-r)Q(u)$ and consider
the change of variables $ u= \frac{u_++u_-}{2} + \frac{u_+-u_-}{2} \sin(\theta)$.  Notice
that $Q(u)\neq0$ on $[ u_-, u_{+}]$.  It follows that
$du=\sqrt{(u- u_{-})(u_{+}-u)}d\theta$ and hence
\eqref{period} can be written in a regularized form as
\[
T(a,E,c)=2\int_{-\frac{\pi}{2}}^{\frac{\pi}{2}}\frac{d\theta}{\sqrt{Q\left(\frac{u_++u_-}{2} + \frac{u_+-u_-}{2} \sin(\theta)\right)}}.
\]
Similarly the mass, momentum, and Hamiltonian of the traveling wave are given
by the first and second moments of this density, i.e.
\begin{align}
M(a,E,c) = \langle u\rangle &= \int_0^T u(x) dx = 2\int_{u_-}^{u_+} \frac{u\; du}{\sqrt{2\left(E-V(u;a,c)\right)}}\nonumber \\
P(a,E,c) = \langle u^2\rangle &= \int_0^T u^2(x) dx = 2\int_{u_-}^{u_+} \frac{u^2\; du}{\sqrt{2\left(E-V(u;a,c)\right)}}\nonumber \\
H(a,E,c) =  \left\langle \frac{u_x^2}{2} - F(u)\right\rangle &=  2\int_{u_-}^{u_+} \frac{E-V(u;a,c) - F(u)}{\sqrt{2\left(E-V(u;a,c)\right)}}\;du\nonumber .
\end{align}
Notice that these integrals are regularized by the same substitution. In
particular one can differentiate the
above expressions with respect to the parameters $(a,E,c)$ and
 the derivatives of these quantities will play an important role
in the subsequent theory.
Note that there is a third constant of integration $x_0$ corresponding to
translation invariance, but this can be modded out and does not play an
important role in the theory.

These quantities satisfy a number of identities, as is derived in the
appendix. In particular if we define the classical action
\begin{equation}
K = \oint p \;dq = \int_0^T u_x^2 \;dx = 2\int_{u_-}^{u_+}\sqrt{2(E-V(u;a,c))}\;du \nonumber
\end{equation}
(which is not itself conserved) then this quantity satisfies the
following relations
\begin{align*}
K_E &= T \\
K_a &= M \\
K_c &= \frac{P}{2}.
\end{align*}
Using the fact that $T,M,P$ and $H$ are $C^1$ functions of parameters $(a,E;c)$,
the above implies the following relationship between the gradients of these
quantities
\[
E \nabla T + a \nabla M + \frac{c}{2} \nabla P + \nabla H = 0
\]
where $\nabla = (\partial_a,\partial_E,\partial_c)$: see the appendix for details of this calculation.
The subsequent theory is developed most naturally in terms of the quantities
$T$, $M$, and $P$.  However, it is possible to restate our results
in terms of $M$, $P$ and $H$ using the above identity.
This is desirable since these have a natural interpretation as conserved
quantities of the partial differential equation.

As noted before this long-wavelength calculation is geometric, and a
number of Jacobian determinants arise. We adopt the following notation
for $2 \times 2$ Jacobians
\[
\{f,g\}_{x,y} = \left|\begin{array}{cc} f_x & g_x \\ f_y & g_y \end{array}\right|
\]
with $\{f,g,h\}_{x,y,z}$ representing the analogous $3\times 3$ Jacobian.

We now begin our study of linear stability of the periodic waves
$u(x)= u(x;a,E,c)$ under small perturbation. To this end, we consider
a small perturbation of the periodic wave $ u(x;a,E)$ of the form
\[
\psi(x,t;a,E,c)= u(x;a,E,c)+\varepsilon v(x,t) +O(\varepsilon^2),
\]
where $0<|\varepsilon|\ll 1$ is a small parameter.  Substituting this
into \eqref{gkdv} and collecting the $O(\varepsilon)$ terms yields
the linearized equation $\partial_x\mathcal{L}[ u]v= -v_t$, where
$\mathcal{L}[u]:=-\partial_{x}^2-f'( u)+c$ is a linear differential operator
with periodic coefficients.  Since the linearized equation is
autonomous in time, we may seek separated solutions of the form
$v(x,t)=e^{-\mu t}v(x)$, which yields the eigenvalue problem
\begin{equation}
\partial_x\mathcal{L}[u]v=\mu v.\label{lineargkdv}
\end{equation}
Note that we consider the linearized operator $\partial_x\mathcal{L}[ u]$ as
a closed linear operator acting on a Banach space $X$ with domain $\mathcal{D}(\partial_x\mathcal{L}[u])$.
In literature, several choices for $X$ have been studied, each of which corresponding to different classes
of admissible perturbations $v$.  In our case, we consider $X=L^2(\RM;\RM)$
and $\mathcal{D}(\partial_x\mathcal{L}[u])=H^3(\RM)$, corresponding to spatially
localized perturbations.
In this case standard Floquet theory
yields the following definitions.

\begin{df}
The monodromy matrix $\MM(\mu)$ is defined to be the period map
\[
\MM(\mu) = \Phi(T,\mu)
\]
where $\Phi(x,\mu)$ satisfies
\begin{equation}
\Phi_x=\HM(x;\mu)\Phi~~~~~~~~~~\Phi(0,\mu)={\bf I}\label{lineargkdvsystem}
\end{equation}
with ${\bf I}$ the $3\times 3$ identity matrix and
\[
\HM(x;\mu) = \left(
             \begin{array}{ccc}
               0 & 1 & 0 \\
               0 & 0 & 1 \\
               -\mu- u_x f''( u) & -f'( u)+c & 0 \\
             \end{array}
           \right).
\]
\end{df}
Given the monodromy the spectrum is characterized as follows:
\begin{df}
We say $\mu\in\spec(\partial_x\mathcal{L}[ u])$ if there exists a
non-trivial bounded function $\psi$ such that $\partial_x\mathcal{L}[
u]\psi=\mu\psi$ or, equivalently if there exists a $\lambda \in \CM$
such that $|\lambda|=1$ and
\[
\det[ \MM(\mu)-\lambda {\bf I} ] =0.
\]
Following Gardner\cite{Gard} we define the periodic Evans function to be
\begin{equation}
D(\mu,\lambda)=\det\left(\MM(\mu)-\lambda {\bf I}\right)\label{evans}.
\end{equation}
Moreover, we say the periodic solution $
u(x;a,E,c)$ is spectrally stable if $\spec(\partial_x\mathcal{L}[ u])$ does
not intersect the open right half plane.
\end{df}

\begin{rem}
Notice that due to the Hamiltonian nature of the problem,
$\spec(\partial_x\mathcal{L}[u])$ is symmetric with respect to
reflections across the real and imaginary axes.  Thus, spectral stability
occurs if and only if $\spec(\partial_x\mathcal{L}[ u])\subset\mathbb{R}i$.
Since we are primarily concerned with roots of $D(\mu,\lambda)$ with $\lambda$
on the unit circle we will frequently work with the function
$D(\mu,e^{i\kappa})$, which is actually the function considered by Gardner.
\end{rem}

In this paper, we will study different asymptotics of this function.  In the next section,
we will study the asymptotics of \eqref{evans} as $\mu\to\infty$.  This will provide
information about the global structure of the spectrum of the linearized
operator $\partial_x\mathcal{L}[u]$, as well as providing us
with a finite wavelength instability index which counts modulo 2 the
number of intersections of the spectrum with the positive real axis.
We then study the asymptotics of
\eqref{evans} in the limit $(\mu,\kappa)\to(0,0)$, which yields
a quantity which we refer to as a modulational stability index, which
is expressed in terms of the derivatives
of the monodromy operator at the origin.

\section{Global Structure of $\spec(\partial_x\mathcal{L}[u])$}

In this section, we review some basic global features of the spectrum
of the linearized operator $\partial_x\mathcal{L}[u]$ which are useful
in a local analysis near $\mu=0$.  We also state
some important properties of the Evans function $D(\mu,\lambda)$ which are vital
to the foregoing analysis.


\begin{prop}\label{globalspec}
The spectrum $\spec\left(\partial_x {\cal L}[u]\right)$ has the following properties:
\begin{itemize}
\item There are no isolated points of the spectrum.  In particular, the spectrum consists of
piecewise smooth arcs.
\item $D(\mu,\lambda)=\det(\MM(\mu)-\lambda {\bf I})=-\lambda^3 + a(\mu) \lambda^2 - a(-\mu)\lambda +1 $ with $a(\mu) = \tr(\MM(\mu)).$
\item The function $a(\mu)$ satisfies $a(0)=3,a'(0)=0.$
\item The entire imaginary axis is contained in the spectrum, i.e. $i \RM \subset \spec(\partial_x\mathcal{L}).$ Further for $|\mu|$ sufficiently large along the imaginary axis the
multiplicity is one.
\item $ \RM \cap \spec(\partial_x\mathcal{L})$ consists of a finite number of
points. In particular there are no bands on the real axis.
\end{itemize}
\end{prop}

\begin{proof}
The first claim, that the spectrum is never discrete, follows from
a basic lemma in the theory of several complex variables: namely
that, if for fixed $\lambda^*$ the function
$D(\mu,\lambda^*)$ has a zero of order $k$ at $\mu^*$ and is holomorphic in a
polydisc about $(\mu^*,\lambda^*)$ then there is some smaller
polydisc about $(\mu^*,\lambda^*)$ so that for every $\lambda$ in a disc about
$\lambda^*$ the function $D(\mu,\lambda)$ (with $\lambda$ fixed) has $k$ roots
in the disc
$|\mu-\mu^*|<\delta$. For details see the text of Gunning\cite{Gun}. It is clear from the implicit function theorem that $\mu$ is a smooth function of
$\lambda$ as long as $\frac{\partial D}{\partial \mu} = \tr(\cof^t(\MM(\mu)-\lambda {\bf I})\MM_\mu^\prime)\neq 0,$ where $\cof$ represents the standard cofactor
matrix.

The second claim is an easy symmetry calculation. The stability problem is
invariant under the
the map $x \mapsto -x, \mu \mapsto -\mu,$ which implies that
\[
\MM(\mu)\sim {\MM}^{-1}(-\mu).
\]
Thus one has
\begin{eqnarray*}
\det[\MM(-\mu)-\lambda {\bf I}]&=&\det[{\MM}^{-1}(\mu)-\lambda {\bf I}]\\
&=&-\lambda^3\det[{\MM}^{-1}(\mu)]\det[\MM(\mu)-\lambda^{-1}]\\
&=&-\lambda^3\left(-\lambda^{-3}+a(\mu)\lambda^{-2}+b(\mu)\lambda^{-1}+1\right)\\
&=&-\lambda^3-b(\mu)\lambda^2-a(\mu)\lambda+1
\end{eqnarray*}
from which it follows $b(\mu)=-a(-\mu)$.

The proof of the third claim will be deferred until lemma \ref{vanishlem}.

The fourth claim is another symmetry argument.
Since $a(\mu)$ is real on the real axis it follows from Schwarz reflection
that for $\mu\in\RM i$, we have
$a(-\mu)=a(\overline{\mu})=\overline{a(\mu)}$ and the characteristic polynomial
takes the form
\[
D(\mu,\lambda)=-\lambda^3+a\lambda^2-\overline{a}\lambda+1
\]
where $a=a(\mu)$, and thus
\[
D(\mu,\lambda)
              = -\lambda^3\overline{D\left(\mu;\overline{\lambda}^{-1}\right)}.
\]
Hence  for imaginary $\mu$ the eigenvalues of the monodromy are symmetric
with respect to the unit circle with the same multiplicities.
Since the monodromy has three eigenvalues, it follows that at least one must lie
on the unit circle.

To see that the multiplicity is eventually one we note that by
 standard asymptotics the monodromy $\MM(\mu)$ satisfies
\[
\MM(\mu) \approx e^{{\bf A}(\mu) T},  ~~~~~~|\mu|\gg 1
\]
where ${\bf A}(\mu)$ is defined by
\[
{\bf A}(\mu)=\left(
    \begin{array}{ccc}
      0 & 1 & 0 \\
      0 & 0 & 1 \\
      -\mu & 0 & 0 \\
    \end{array}
  \right).
\]
The three eigenvalues of $e^{{\bf A}(\mu)T}$ are given by
\begin{equation}
\lambda_1=e^{-\mu^{1/3}T},\;\;\lambda_2=e^{-\mu^{1/3}\omega
T},\;\;\rm{and}\;\;\lambda_3=e^{-\mu^{1/3}\overline{\omega}T}\label{assymev}
\end{equation}
where $\omega = e^{2\pi i/3}$ is the principle third root of
unity.  If $\mu\in\RM^+ i$ it follows that
$\lambda_1=\exp\left(-|\mu|^{1/3}e^{i\pi /6}T\right)$
and since $\cos(\pi/6)>0$ we have $|\lambda_1|\to 0$ as $\RM^+
i\ni\mu\to\infty$.  Similarly,
$\lambda_2=\exp\left(-|\mu|^{1/3}e^{5\pi/6}T\right)$ and
$\lambda_3=\exp\left(|\mu|^{1/3}i\right)$ so that $|\lambda_2|\to
\infty$ as $\RM^+ i\ni\mu\to\infty$ and $|\lambda_3|=1$.  Thus, for
$\mu\in\RM^+i$ large, we have that $\mu$ is an eigenvalue of
multiplicity one. Similarly, we can show $|\lambda_1|\to \infty$,
$|\lambda_3|\to 0$ as $\RM^+i\ni\mu\to-\infty$ and
$|\lambda_2|=1$ for $\mu\in\RM^-i$, $|\mu|\gg 1$.  Therefore, it
follows that $\mu\in\spec(\partial_x\mathcal{L}[ u])$ with multiplicity one for
$\mu\in\RM i$, $|\mu|\gg 1$.

The final claim follows from a similar asymptotic calculation together
with an analyticity argument. Notice that for $\mu$ real the eigenvalues of
the monodromy are either
all real or one real and one complex conjugate pair. If the eigenvalues
lie on the unit circle then in the first case $1$ or $-1$ must be an
eigenvalue. In the second one must have a complex conjugate pair of
eigenvalues, and thus (since the determinant is one) $1$ must be an eigenvalue.
Thus if a point on the real axis is in the spectrum then either
$\det(\MM(\mu)-{\bf I})$ or $\det(\MM(\mu)+{\bf I})$ must vanish. Since
$\det(\MM(\mu)\pm {\bf I})$ are entire
functions it follows that either they are identically zero or the
zero set has not finite accumulation points.
The large $\mu$ asymptotics implies that they cannot be identically zero,
therefore the zero set must be discrete.
Further the large $\mu$ asymptotics implies that for sufficiently large
$\mu$ along the real axis $\mu \notin \spec\left(\partial_x \mathcal{L}[u]\right)$,
so the spectrum is confined to a compact subset of the real line,
and there are only a finite number of real eigenvalues.
\end{proof}

\begin{rem}
Note that, in the calculation of H\v{a}r\v{a}gu\c{s} and Kapitula\cite{HK} the real eigenvalues play a slightly different role than other eigenvalues off of the
imaginary axis. The fact that there are only a finite number of these indicates
that there are only a finite number of values of the Floquet parameter for which there are real eigenvalues: $\kappa_r(\gamma)=0$ for all but a finite
number of values of the Floquet parameter $\gamma$ in their notation.
\end{rem}

\subsection{Analysis of $\spec(\partial_x\mathcal{L}[u])\bigcap\RM$}

We now move on to study the structure of $\spec(\partial_x\mathcal{L}[
u])\cap\RM$ more carefully.  Suppose that $\mu\in\RM$.
Clearly $1$ being an eigenvalue of $\MM(\mu)$ is a sufficient
condition for
$\mu\in\spec(\partial_x\mathcal{L}[ u])$ and
thus vanishing of $D(\mu,1)=a(\mu)-a(-\mu)$ is a sufficient condition
$\mu\in\spec(\partial_x\mathcal{L}[ u])\cap\RM$.  Notice that by the translation
invariance of \eqref{lineargkdv} we have $D(0,1)=0$ by Noether's theorem.
The question is whether $D(\mu,1)$ has any
other real roots.  If it does, then the eigenvalue problem
\eqref{lineargkdv} is spectrally unstable, due to the presence of
a real non-zero element of $\spec(\partial_x\mathcal{L}[ u])$.  In order to detect
this instability, we calculate the orientation index
\[
D(\infty,1)D_{\mu\mu\mu}(0,1).
\]
As we will show, the negativity of this index is sufficient to imply a non-trivial
intersection of $\spec(\partial_x\mathcal{L}[u])$ with the real line.

\begin{lem}
The function $D(\;\cdot\;,1):\RM\to\RM$ is an odd function which satisfies the
asymptotic relation
\[
\lim_{\RM\ni\mu\to\pm\infty} D(\mu,1)=\mp\infty.
\]
\end{lem}

\begin{proof}
Clearly, $D(\;\cdot\;,1)$ is an odd function of its argument, and hence it is
sufficient to consider the limit as $\mu\to\infty$.  To begin, define a new variable
$\rho=\mu^{1/3}T$.  Then from the asymptotic relations \eqref{assymev} we have
\begin{eqnarray*}
a(\rho)&=&e^{-\rho}+e^{-(-1+\sqrt{3}i)\rho/2}+e^{-(-1-\sqrt{3}i)\rho/2}\\
\widetilde{a}(\rho)&=&e^{\rho}+e^{-(1+\sqrt{3}i)\rho/2}+e^{-(1-\sqrt{3}i)\rho/2}
\end{eqnarray*}
where $\widetilde{a}(\rho)$ is the trace when you take $\mu\to-\mu$.  It follows that
$D(\mu,1)=a(\rho)-\widetilde{a}(\rho)$ behaves like $-e^{\rho}$ for large positive
$\rho$, i.e. $\mu\gg 0$.  This completes the proof.
\end{proof}

From these results, we have the following theorem relating the
sign of $\tr(\MM_{\mu\mu\mu}(0))$ to the stability of the underlying
periodic wave.

\begin{thm}\label{realspectrum}
If $a'''(0)=\tr(\MM_{\mu\mu\mu}(0))>0$, then the number of roots of $D(\mu,1)$
(i.e. the number of periodic eigenvalues)
on the positive real axis is odd. In particular
$\spec(\partial_x\mathcal{L}[ u])\cap\RM^*\neq\emptyset$ and
the eigenvalue problem
\eqref{lineargkdv} is spectrally unstable.
\end{thm}

\begin{proof}
We show in Lemma 2 that $D(0,1)=D_\mu(0,1)=D_{\mu\mu}(0,1)=0$ and
$D_{\mu\mu\mu}(0,1)=2a'''(0)$.  Thus, if $a'''(0)>0$, then 
$D(\mu;1)$ is positive for small positive values of $\mu$.
Since $D(\mu,1)$ is negative for sufficiently large $\mu$
we know that  $D(\pm\mu^*,0)=0$ for some $\mu^*\in\RM\setminus\{0\}$, which completes the proof. In the next section we establish the
following formula for $D_{\mu\mu\mu}(0,1)$, the first non-vanishing derivative
\[
D_{\mu\mu\mu}(0,1) = -3\left| \begin{array}{ccc} T_a & M_a & P_a \\
 T_E & M_E & P_E \\
 T_c & M_c & P_c  \end{array}\right|
 =6 \left| \begin{array}{ccc} K_{aa} & K_{aE} & K_{ac} \\
 K_{aE} & K_{EE} & K_{cE} \\
 K_{ac} & K_{cE} & K_{cc}  \end{array}\right|
  = \frac{3}{E}
\left| \begin{array}{ccc} M_a & P_a & H_a \\
 M_E & P_E & H_E \\
 M_c & P_c & H_c  \end{array}\right|
 \]
where again $K$ is the classical action of the traveling wave ODE.
Hence this ``orientation index'' can be expressed in terms of the
Jacobian of the map between the constants of integration of the
traveling wave ordinary differential equation $(a,E;c)$ and the
conserved quantities of the g-KdV $(M,P,H).$ This orientation
index is analogous to the quantity which is calculated in the stability
theory of the solitary waves.
\end{proof}

It is important to notice the instability detected by
Theorem \ref{realspectrum} is an instability with respect to
\textit{finite (bounded) wavelength} perturbations.  In the next section
we will derive a \textit{modulational} stability index which detects
instability with respect to arbitrarily long wavelength perturbations.
See the comments at the end
of the section 3. The solitary wave solutions go unstable in the manner
detected by Theorem \ref{realspectrum}, through
the creation of a pair of eigenvalues on the real axis. In general the periodic
waves seem to first go unstable through the creation of a curve of spectrum
which does not  intersect the real axis, and later  there
is a secondary bifurcation resulting in a real eigenvalue. This phenomenon
appears to have first been observed by Kapitula and H{\v a}r{\v a}gu{\c s},
who established that small amplitude periodic waves first go unstable at
$p=2$, as compared with $p=4$ for the solitary waves.
While we don't have a general proof of this we do show that, in the case of power
law nonlinearity, there is a real periodic eigenvalue as well as a
band of unstable spectrum connected to the origin. It is also worth
noting that the analogous calculation for $D(\mu,-1)$ shows that the
number of anti-periodic eigenvalues on the real axis is always even. While
this is not
useful for proving the existence of instabilities it does eliminate some
possible modes of instability.

\section{Local Analysis of the Period Map}

In this section, we turn our attention to studying the monodromy
map $\MM(\mu)$ near the origin.  To this end, we determine
the asymptotic behavior of $D(\mu,e^{i\kappa})$ as $\mu\to 0$.  We begin by proving
that $D(0,e^{i\kappa})$ has a zero of multiplicity 
three at $\kappa=0$. It follows by directly computing the Jordan normal form of $\MM(0)$
that $\lambda=1$ is an eigenvalue of algebraic  multiplicity three and
geometric multiplicity (generically) two.
This fact reflects the following structure in the manifold of solutions to
the ordinary differential equation defining the traveling waves: the traveling waves form a three parameter manifold, with traveling waves of constant period forming a two parameter submanifold. The two eigenvectors of the period map
correspond to elements of the tangent plane to the submanifold of constant
period solutions,
while the third vector in the Jordan chain is associated to the normal
to the constant period submanifold.

Using perturbation theory appropriate to a Jordan block, as well
as the Hamiltonian symmetry inherent in \eqref{lineargkdv}, we prove
the three roots of $D(\mu,e^{i\kappa})$ bifurcate from $\mu=0$ analytically in
$\kappa$
in a neighborhood of $\kappa=0$, and derive
a necessary and sufficient condition for modulational instability of the underlying
periodic wave $u(x;a,E,c)$ in terms of derivatives of the monodromy operator.
Note that this conclusion is somewhat unexpected: normally the eigenvalues
of a non-trivial Jordan block do not bifurcate analytically but instead
admit a Puiseaux series in fractional powers. However because of the
symmetries of the
problem the admissible perturbations are severely restricted, resulting in
a non-generic bifurcation.

\subsection{Calculation of the Period Map}

The first major calculation we present is an explicit calculation of
the monodromy matrix at the origin in terms of the derivatives of
the underlying periodic solution $u$ with respect to the parameters.
We do this by first computing a matrix valued solution to the
ordinary differential equation satisfying the wrong initial
condition: ${\bf U}(0,0)$ is non-singular but not the identity.
One can then multiply on the right by ${\bf U}^{-1}(0,0)$ to find the
monodromy matrix. We find that (as expected) the monodromy operator
$\MM(\mu)$ has a non-trivial Jordan form when $\mu=0$.  Our goal is then
to utilize perturbation theory of Jordan blocks to calculate the normal
form of the characteristic polynomial in a neighborhood of
$\mu=0,\;\lambda=1$, where $\lambda$ is the eigenvalue parameter
of the monodromy operator.

To begin we write the above
third order eigenvalue problem as a first order system as in \eqref{lineargkdvsystem}. In particular, notice that $\tr(\HM(x;\mu))=0$ for all $x,\mu$, and
thus $\det(\Phi(x;\mu))=1$ for all
$\mu\in\mathbb{C}$, implying $\det(\MM(\mu))=1$. In order to calculate a matrix solution
$\Phi(x;\mu)$, we must first find three linearly independent solutions of the
above system. In general, this is a daunting task, but since
the above system with $\mu=0$ arises as the Frechet derivative (linearization)
of an integrable ordinary differential equation this can be done by
considering infinitesimal variations of the constants of integration
in the defining ordinary differential equation, and thus generating
the tangent space. As noted earlier
the solutions $ u(x-x_0;a,E,c)$ constitute a 4-dimensional solution
manifold of \eqref{gkdv} parameterized by
$x_0,a,E,c$. The solutions  of the linearized
operator space are given by the generators $\frac{d}{dx}$,
$\frac{d}{da}$, and $\frac{d}{dE}$ acting on the solution $u(x;a,E,c)$.
The action of  the generator $\frac{d}{dc}$ is somewhat different and is
connected with the generalized null-space. This will become important in the next section.

\begin{prop}\label{Kernel}Let $u(x;a,E,c)$ be the solution to the traveling
wave equation \eqref{dphi2} satisfying $u(0;a,E,c)=u_-,u_x(0;a,E,c)=0.$
A basis of solutions to the third order system
\[
Y_x = \HM(x;0) Y
\]
is given by
\begin{align*}
&Y_1^t = (u_x,u_{xx},u_{xxx})  \\
&Y_2^t = (u_a,u_{ax},u_{axx})  \\
&Y_3^t = (u_E,u_{Ex},u_{Exx}).
\end{align*}
A particular solution to the inhomogeneous problem
\[
Y_x = \HM(x;0) Y + W
\]
 where $W^t = (0,0,u_x)$ is  given by
\[
Y_3^t = (u_c,u_{cx},u_{cxx}).
\]
\end{prop}
\begin{proof}
A straightforward calculation.  Notice that it follows that $\partial_x\mathcal{L}[u]\left(-u_{c}\right)=u_x$.
\end{proof}
The fact that $u_a,u_E$ are not periodic - they exhibit secular growth due to
the variation of the period with respect to the parameters - gives an
indication that the eigenspaces of the monodromy at $\mu=0$ are not
semi-simple, and hence we expect the existence of a non-trivial Jordan block of the monodromy
map $\MM(0)$.

By the above proposition, three linearly independent solutions of
\eqref{lineargkdvsystem} corresponding to $\mu=0$ are given by
\begin{equation}
Y_1(x)= \left(
  \begin{array}{c}
     u'(x;a,E,c) \\
     u''(x;a,E,c) \\
     u'''(x;a,E,c) \\
  \end{array}
\right) \;\;\; Y_2(x) = \left(
           \begin{array}{c}
              u_a(x;a,E,c) \\
              u_a'(x;a,E,c) \\
              u_a''(x;a,E,c) \\
           \end{array}
         \right)
\;\;\;Y_3(x) = \left(
  \begin{array}{c}
     u_E(x;a,E,c) \\
     u_E'(x;a,E,c) \\
     u_E''(x;a,E,c) \\
  \end{array}
\right).\label{solutionvec}
\end{equation}
By hypothesis, for any $a,E,c\in\mathbb{R}$ the solution $ u$
satisfies
\begin{eqnarray}
 u(0;a,E,c) &=&  u_{-} \;=\;  u(T;a,E,c)\label{initial1}\\
\partial_x u(0;a,E,c) &=& 0 \;=\; \partial_x u(T;a,E,c)\label{initial2}\\
\partial_{xx} u(0;a,E,c) &=& a-f(u_{-})+cu_- \;=\; \partial_{xx} u(T;a,E,c)\label{initial3}.
\end{eqnarray}
Moreover, from equation \eqref{gkdv} it follows that
\[
u_{xxx}(0;a,E,c)=cu_x(0;a,E,c)-\frac{d}{dx}\left(f( u(x;a,E,c))\right)\big{|}_{x=0}
=0.
\]
Defining $\UM(x,0)=[Y_1(x),Y_2(x),Y_3(x)]$ to be the corresponding solution matrix, then direct calculations yield
\begin{equation}
\UM(0,0)=\left(\begin{array}{ccc}
0 & \partial_a u_{-} & \partial_E u_{-}\\
a-f(u_{-})+cu_- & 0 & 0\\
0 & 1+(c-f'(u_{-}))\partial_a u_{-} & (c-f'( u_{-}))\partial_E u_{-}
\end{array}\right).\label{MatrixSolutionInitial}
\end{equation}
Note that differentiating the relation $E-V(u_-)=0$ gives  the relation
$-V'(u_{-})\partial_Eu_-=\det(\UM(0,0))=-1$, so these solutions are linearly
independent at $x=0$, and hence for all $x$. Thus we can compute ${\bf U}(T,0)$
and right-multiply by ${\bf U}^{-1}(0,0)$ to give the monodromy ${\bf M}(0)$

The matrix $\UM(T,0)$ can be calculated by differentiating
\eqref{initial1}-\eqref{initial3}
with respect to the parameters $a$ and $E$ by use of the chain rule.
For example, differentiating the relation \eqref{initial1} with respect to the parameter
$E$ gives
\[
\partial_E  u(T;a,c,E)+\frac{\partial u}{\partial x}\left(T;a,E,c\right) T_E(a,c,E) = \partial_E u_{-}.
\]
Since the derivative vanishes at the period points this implies $\partial_E u(T)=\partial_E u_{-}$.  Continuing in this manner gives the following expression
for the change in tis matrix solution across the period:
\begin{equation}
\UM(T,0) = \UM(0,0)+ \left(\begin{array}{ccc}
0 & 0 & 0\\
0 & V'(u_{-};a,c)T_a & V'(u_{-};a,c)T_E\\
0 & 0 & 0
\end{array}\right).\label{MatrixSolutionDifference}
\end{equation}
In particular, we find that $\UM(T,0)-\UM(0,0)$ is a rank one matrix,
which naturally leads to the following proposition.

\begin{prop}\label{monodromy}
There exists a basis in $\mathbb{R}^3$ such that the monodromy matrix
$M(\mu)$ evaluated at
$\mu=0$ takes the following Jordan normal form:
\begin{equation}
\MM(0) \sim \left(\begin{array}{ccc}
1 & 0 & 0\\
0 & 1 & \sigma\\
0 & 0 & 1
\end{array}\right)\label{jordanform}
\end{equation}
where $\sigma\neq 0$ as long as $T_a$ and $T_E$ do not simultaneously
vanish. In particular,
the monodromy operator at $\mu=0$ has a single eigenvalue of
$\lambda=1$ with algebraic multiplicity three and geometric
multiplicity two as long as the period is not at a critical point with respect to the parameters $a,E$ for fixed wavespeed $c$.
\end{prop}

\begin{proof}
Recall $\det(\UM(0,0))=-V'(u_{-})\partial_E u_{-}=-1$,
so $\UM(0,0)$ is invertible. Multiplying the above expression
on the right by the matrix $\UM^{-1}(0,0)$ yields the monodromy
matrix at the origin
\[
\MM(\mu=0):=I+\vec{w}\otimes\vec{v}\;\UM^{-1}(0)
\]
where $\vec{w}=(0,1,0)^{T}$ and $\vec{v}=\left(0,V'(u_{-})T_a,V'(u_{-})T_E\right)^{T}$.
Next, notice that
\[
\UM(0,0)\left(\begin{array}{ccc}
0 & -T_a & -T_E\\
0 & 0 & 0\\
0 & 0 & 0
\end{array}\right) = \vec{w}\otimes\vec{v}
\]
and hence defining $\NM:=\UM^{-1}(0)\MM(0)\UM(0)$ gives the equation
\[
\NM=\left(\begin{array}{ccc}
1 & -T_a & -T_E\\
0 & 1 & 0\\
0 & 0 & 1
\end{array}\right).
\]
It follows that
\[
Ker(\NM-{\bf I})=span\left\{(1,0,0)^T,\;(0,T_E,-T_a)^T\right\}
\]
Now, take $\vec{v_3}:=(0,-T_a,-T_E)$ and notice that $v_3\notin
Ker(\NM-{\bf I})$.  The Jordan structure then follows by noticing then
that $(\NM-{\bf I})\vec{v_3}=(T_a^2 + T_E^2)(1,0,0)^T \in \ker({\bf N}-{\bf I})$.

For traveling waves below the separatrix a result of Schaaf (see appendix)
shows that\footnote{Under some mild assumptions on the nonlinearity, which are satisfied for the power law nonlinearity.}
 the period is a strictly increasing function of the energy, $T_E>0$
and thus the genericity condition is always met. In other situations we will
assume that this condition is met unless otherwise stated.
\end{proof}

\subsection{Asymptotic Analysis of $D(\mu,\kappa)$ near $(\mu,\kappa)=(0,0)$}

We now analyze the characteristic polynomial of $\MM(\mu)$ in a neighborhood
of $\mu=0$ by considering $\MM(\mu)$ as a small perturbation of the matrix
$\MM(0)$ constructed above.  It is well understood how the eigenvalues of
a Jordan block bifurcate under perturbation: see Kato\cite{K} or Moro, Burke and Overton\cite{MBO}. It is worth noting, however, that in this case the
bifurcation is highly non-generic due to the constraints imposed by the
symmetry of the problem.

Recall from Proposition \ref{globalspec} that the spectrum near $\mu=0$
is continuous. By the analyticity of $\MM(\mu)$ in a neighborhood of
$\mu=0$, we can expand $\MM(\mu)$ for small $\mu$ as
\[
\MM(\mu)=\MM(0)+\mu
\MM_\mu(0)+\frac{\mu^2}{2}\MM_{\mu\mu}(0)+\mathcal{O}\left(|\mu|^3\right)
\]
where $\MM_\mu(0)=[M^{(1)}_{i,j}]$ and $\MM_{\mu\mu}(0)=[M^{(2)}_{i,j}]$.  If one
makes a similarity transform $\widetilde\MM(\mu) = {\bf V}^{-1}\MM(\mu) {\bf V}$
so that
$\widetilde \MM(0)$ is in the Jordan normal form \eqref{jordanform} then a
direct calculation using the above second order expansion of
$\widetilde \MM(\mu)$ implies that in a neighborhood of $\mu=0$,  the
characteristic polynomial can be expressed as
\begin{eqnarray}
D(\mu,e^{i\kappa})&=&\det\left((\widetilde \MM(\mu)-I)-(e^{i\kappa}-1){\bf I}\right)\nonumber\\
&=&-\eta^3+\eta^2\left(\mu \tr\left(\widetilde\MM_\mu(0)\right)+\frac{\mu^2}{2}\tr(\widetilde\MM_{\mu\mu}(0))\right)\nonumber\\
&-&\eta\left(\mu {\widetilde M}^{(1)}_{3,2}\sigma+\mu^2\left(\frac{1}{2}\left(\tr(\widetilde\MM_\mu)\right)^2 -\frac{1}{2} \tr(\widetilde\MM_\mu^2)-\frac{\sigma}{2}{\widetilde M}^{(2)}_{3,2}\right)\right)\nonumber\\
&-&\sigma(\widetilde M^{(1)}_{1,1} \widetilde M^{(1)}_{3,2}-\widetilde M^{(1)}_{3,1}\widetilde M^{(1)}_{1,2})\mu^2\nonumber\\
&+&\mu^3\left(\det\left(\widetilde\MM_\mu(0)\right)+\sigma S\right)
+\mathcal{O}\left(4\right),\label{perturbation1}
\end{eqnarray}
where $\eta=e^{i\kappa}-1$, $S$ represents mixed terms from $\widetilde\MM_\mu(0)$ and
$\widetilde\MM_{\mu\mu}(0)$, $\sigma$ is as in Proposition \ref{monodromy}, and the notation $\mathcal{O}\left(4\right)$ represents
terms whose degree is four or higher. Notice there are no other $\mu^3$ terms since $\MM(0)-I$ has rank one.
Our next goal is to determine the dominant balance of the equation $D(\mu,e^{i\kappa})=0$
in a neighborhood of $(\mu,\kappa)=(0,0)$.

A useful construction for implicit function calculations of this type is that
of the Newton diagram, which is a subset of the non-negative integer lattice.
 A vertex $(i,j)$ is included if the coefficient of $\eta^{3-i} \mu^j$
in \eqref{perturbation1} is {\em non-zero}, otherwise the vertex is not
included. The lower convex hull of the Newton  diagram is made up of a
collection of line segments.
For each line segment of the lower convex hull let $m$ be the
horizontal length of the segment and $s$ the slope of the segment.
Corresponding to each such line segment there are $m$
 distinct solution branches of the form
\[
\eta_k(\mu) = \sum_{i} \alpha^{(k)}_i \mu^{s i},
\]
where $\alpha^{(k)}_1\neq 0$ and $k$ ranges from $1$ to $m$.
For details see the book of Baumgartel\cite{Baum} or Hilton\cite{Hilton}. This is equivalent to the method of
``dominant balance'' presented in
textbooks on asymptotic methods  but is somewhat more systematic.
For instance, in our case if the coefficient of the $\eta^1\mu^1$
term (-$\sigma \tilde M^{(1)}_{3,2}$) is non-vanishing then there are two solution
branches in which $\eta$ has an expansion in powers of $\mu^{\frac{1}{2}}$ and
one with an expansion in integer powers. These correspond to the breaking
up of the $2\times 2$ and $1\times 1$ Jordan blocks respectively.
However as mentioned above the symmetry $\MM(\mu) \sim \MM^{-1}(-\mu)$ causes
a number of terms in \eqref{perturbation1} to vanish, which
leads to an expansion in integer powers of $\mu$. This is the content of the
next lemma.

\begin{lem}
\label{vanishlem}
The equation $D(\mu,e^{i\kappa})=0$ has the
following normal form in a neighborhood of $(\mu,\kappa)=(0,0)$:
\begin{eqnarray*}
-(i\kappa)^3+\frac{i\kappa\mu^2}{2}\tr\left(\MM_{\mu\mu}(0)\right)+\frac{\mu^3}{3}\tr\left(\MM_{\mu\mu\mu}(0)\right) + O(4)=0
\end{eqnarray*}
whose Newton diagram is depicted in Figure 2.
\end{lem}

\begin{figure}
\centering
\includegraphics[scale=.35]{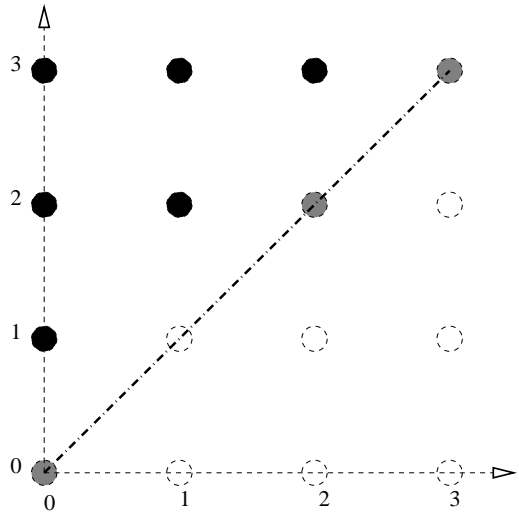}
\caption{The Newton diagram corresponding to the asymptotic
expansion of $D(\mu,e^{i\kappa})=0$ in a neighborhood of
$(\mu,\kappa)=(0,0)$ is shown to $O(|\mu|^3)$.  Terms associated
to open circles are shown to vanish due to the natural symmetries
inherent in \eqref{gkdv}.  The grey circles are non-vanishing terms which
are a part of the lower convex hull. The black circles lie above the
lower convex hull and thus do not contribute to the leading order
asymptotics.}
\end{figure}

\begin{proof}
Define functions a, b, c on a neighborhood of $\mu=0$ by
\begin{equation}
\det[(\MM(\mu)-I)-(e^{i\kappa}-1)
I]=-\eta^3+(a(\mu)-3)\eta^2+b(\mu)\eta+c(\mu).\label{charpoly}
\end{equation}
where $\eta=e^{i\kappa}-1$.  Notice in particular that $\eta=i\kappa+\mathcal{O}\left(\kappa^2\right)$
in a neighborhood of $\kappa=0$.  By \eqref{perturbation1}, it follows
\begin{eqnarray*}
a(\mu)&=& \tr(\MM(\mu))=3+\mu \tr(\MM_{\mu}(0))+\frac{\mu^2}{2}\tr(\MM_{\mu\mu}(0))+\frac{\mu^3}{6}\tr(\MM_{\mu\mu\mu}(0))+\mathcal{O}(|\mu|^4)\\
b(\mu)&=& \frac{1}{2}\left(\tr((\MM(\mu)-I)^2)-\tr(\MM(\mu)-I)^2\right)=-\mu M^{(1)}_{3,2}\sigma-\mu^2\left(\frac{1}{2} \tr(\MM_\mu)^2 -\frac{1}{2} \tr(\MM_\mu^2) -\frac{\sigma}{2}\tilde M_{3,2}\right)+\mathcal{O}(|\mu|^3)\\
c(\mu)&=& \det(\MM(\mu)-I)=-\sigma
(M^{(1)}_{1,1}M^{(1)}_{3,2}-M^{(1)}_{3,1}M^{(1)}_{1,2})\mu^2+\left(\det(\MM_{\mu}(0))+\sigma
S\right)+\mathcal{O}(|\mu|^4)
\end{eqnarray*}
Using the symmetry $\MM(-\mu)\sim \MM(\mu)^{-1}$, we have
\begin{eqnarray*}
c(-\mu)&=&\det[\MM(-\mu)-I]\\
&=&\det[\MM(\mu)]^{-1}\det[I-\MM(\mu)]\\
&=&-\det[\MM(\mu)-I]\\
&=&-c(\mu)
\end{eqnarray*}
since $\det[\MM(\mu)]=1$ for all $\mu\in\mathbb{C}$.  Hence $c$ is an
odd function of $\mu$. Also, since $\MM(0)-{\bf I}$ has rank one,
\eqref{perturbation1} along with the above analysis implies that
$c(\mu)=\mathcal{O}(|\mu|^3)$, from which it follows
$c''(0)=-2\sigma(M^{(1)}_{1,1}M^{(1)}_{3,2}-M^{(1)}_{3,1}M^{(1)}_{1,2})=0$.

Similarly, using \eqref{charpoly} we have
\begin{eqnarray*}
\det[\MM(\mu)-\lambda {\bf I}]&=&-\lambda^3 \det\left[\MM(-\mu)-\frac{1}{\lambda}\right]\\
&=& -\lambda^3\left(\left(\frac{1}{\lambda}-1\right)^3 +
a(-\mu)\left(\frac{1}{\lambda}-1\right)^2
+b(-\mu)\left(\frac{1}{\lambda}-1\right) + c(-\mu)\right)\\
&=&-(\lambda-1)^3 - a(-\mu)\lambda
(\lambda-1)^2+b(-\mu)\lambda^2(\lambda-1)- c(-\mu)\lambda^3.
\end{eqnarray*}
Comparing the $\lambda^2$ and $\lambda^3$ terms above with those
in \eqref{charpoly} we get the relations
\[
\left\{
  \begin{array}{ll}
    b(\mu)=2a(\mu)-a(-\mu)-3, & \hbox{and} \\
    a(\mu)-b(\mu)+c(\mu)=3. &
  \end{array}
\right.
\]
Since $c(\mu)=\mathcal{O}(|\mu|^3)$, these relations imply
$\tr\left(\MM_{\mu}(0)\right)=a'(0)=0$ and $-\sigma \tilde M^{(1)}_{3,2}=b'(0)=0$.
By recalling $\sigma\neq 0$ from \ref{monodromy}, this
implies $\tilde M^{(1)}_{3,2}=0$. Moreover, we know that $a''(0)=b''(0)$ and
hence $b''(0)=\tr(\MM_{\mu\mu}(0))$. Also, we have $b'''(0)=3a'''(0)$
and $c'''(0)=b'''(0)-a'''(0)=2a'''(0)$ and hence
$c'''(0)=2\tr(\MM_{\mu\mu\mu}(0))$. The corollary follows by
analyzing equation \eqref{perturbation1} as well as the
corresponding Newton diagram (see Figure 1).
\end{proof}

From this it follows that, in the neighborhood of the origin, the
leading order piece of the periodic Evans function is a
homogeneous cubic polynomial in $\kappa, \mu$. The implicit function
theorem fails, but in a trivial way that is easily corrected, leading to
the following theorem:

\begin{figure}
\begin{center}
\includegraphics[scale=.35]{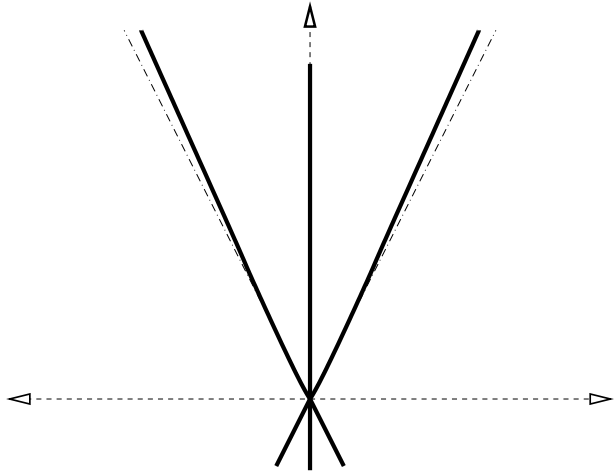}
\end{center}
 \caption{When $\Delta(f; u)<0$, the local normal form of
$\spec(\partial_x \mathcal{L}[ u])$ consists of a segment of the imaginary
axis union with two straight lines making equal angles with the
imaginary axis.  Notice that these lines intersect at the origin,
corresponding to the fact that $1$ is an eigenvalue of $\MM(0)$ with algebraic
multiplicity three.}
\end{figure}

\begin{thm}\label{discriminantthm}
With the above notation, define
\begin{equation}
\Delta(f;
u)=\frac{1}{2}\left(\tr(\MM_{\mu\mu}(0))\right)^3-3\left(\tr(\MM_{\mu\mu\mu}(0))\right)^2,\label{discriminant}
\end{equation}
where $f$ denotes the dependence on the non-linearity used in
\eqref{gkdv}, and suppose $\tr(\MM_{\mu\mu\mu}(0))\neq 0.$ If $\Delta>0$,
then the imaginary axis in the
neighborhood in the origin is in the spectrum with multiplicity
three. If $\Delta<0$ then the imaginary axis in a neighborhood
of the origin is in the spectrum with multiplicity one, together
with two curves which are tangent to lines through the origin with
angle $\arg(i y_{2,3})$ - see Figure 2. In particular in the latter case
the periodic wave is modulationally unstable.
\end{thm}

\begin{proof}
Since the leading order piece of the Evan's function is homogeneous it
suggests working with a projective coordinate $y =\frac{i \mu}{\kappa}.$
Making such a change of variables leads to the equation
\[
1 - \frac{y^2}{2}\tr\left(\MM_{\mu\mu}(0)\right)+\frac{y^3}{3}\tr\left(\MM_{\mu\mu\mu}(0)\right) + \kappa E(\kappa,y) = 0
\]
where $E(\kappa,y)$ is continuous in a neighborhood of the origin.
The implicit function theorem applies in a neighborhood of
$(y=y_{1,2,3},\kappa=0)$ as long as the roots $y_{1,2,3}$ of the
above cubic in $y$ are distinct, which is true as long as the
discriminant $\Delta$ is not zero. In terms of the original
variable $\mu$ we have the three solution branches
\[
\mu_{1,2,3} = i y_{1,2,3} \kappa + O(\kappa^2)
\]
This cubic has three real roots
when $\Delta>0$, giving three branches of spectrum emerging from the
origin tangent to the imaginary axis. It is clear from symmetry that
these must in fact lie on the imaginary axis, giving a interval of spectrum
of multiplicity three along the imaginary axis. In the case that the
discriminant is negative there is one real root and two complex
conjugate roots,
giving one branch of spectrum along the imaginary axis and two branches
emerging in the complex plane.
\end{proof}

\begin{rem}
First we note that $\tr(\MM_{\mu\mu}(0))<0$ is a sufficient condition for modulational
instability of the periodic wave.

Secondly we note that the Newton diagram is independent of
$\tr(\MM_{\mu\mu}(0))$ but the dominant balance changes if
$\tr(\MM_{\mu\mu\mu}(0))$ should happen to vanish. Later we will
give a formula for $\tr(\MM_{\mu\mu\mu}(0))$ in terms of the Jacobian
of a map, and we will see that vanishing of this quantity signals
a change in the Jordan structure of the underlying linearized operator.

Finally we note that this agrees with the result of Bottman and
Deconinck\cite{BD}, in which they considered cnoidal wave solutions
to the KdV.  Using the algebro-geometric techniques of Belokolos, Bobenko,
Enol'skii, Its and Matveev\cite{BBEIM} they explicitly computed the
spectrum of the linearized operator and found that such solutions are
always linearly stable. Their results prove that
an interval of the imaginary axis containing the origin is
in the spectrum of the linearized operator, with multiplicity
three.  Our results imply this is a generic phenomenon: either
one has an interval of spectrum of multiplicity three about  the
origin, or one has three curves intersecting at the origin.
In the  KdV case the discriminant $\Delta$ is expressible in terms of
elliptic functions in this case and must be positive,
although we have not tried to show this.
\end{rem}

It is instructive to compare this theorem with that of Theorem
\ref{realspectrum}.  By Theorem \ref{discriminantthm}, the sign
of $\tr(\MM_{\mu\mu\mu}(0))$ has no effect on the spectral stability
of the underlying periodic wave in a sufficiently small
neighborhood of the origin.  However, Theorem \ref{realspectrum}
guarantees the existence of unstable real spectrum if
$\sgn\left(\tr(\MM_{\mu\mu\mu}(0))\right)>0$.  To reconcile these
results, notice Proposition \eqref{globalspec} implies there is no
unstable real spectrum sufficiently close to the origin.  Thus,
the instability brought on by $\tr(\MM_{\mu\mu\mu}(0))>0$ is not
local to $\mu=0$, and hence should not be detected  by the
quantity $\Delta(a,E,c)$.

Our next goal is to try to evaluate the modulational stability index
$\Delta(a,E,c)$ as well as the finite wavelength orientation index $\tr(\MM_{\mu\mu\mu}(0))$
in terms of the conserved quantities of \eqref{gkdv}.
This can be done very explicitly.  Notice that while we have chosen to express the
coefficients as $\tr(\MM_{\mu\mu}(0))$ and $\tr(\MM_{\mu\mu\mu}(0))$, which
suggests that they arise at second order and third order in a perturbation
calculation for small $\mu$, these quantities can be expressed in
terms of quantities which arise at  first and second order in $\mu$
due to the invariance of the problem under the map $\mu \rightarrow -\mu$.
Further while all of the first order terms contribute only a few
terms which are second order actually contribute - these are
the terms which are associated to the minors of the
off-diagonal piece of the unperturbed Jordan form. These second order terms are
explicitly computable via a single quadrature.

\begin{thm}\label{TraceIdentities}
We have the following identities:
\begin{align}
\tr(\MM_{\mu\mu})|_{\mu=0} &= \{T,P\}_{E,c} +2 \{M,P\}_{a,E} \nonumber \\
\tr(\MM_{\mu\mu\mu})|_{\mu=0}&= -\frac{3}{2}\{T,M,P\}_{a,E,c}\nonumber
\end{align}
where $T,M,P$ are the period, mass, and momentum of the underlying traveling
wave and $a,E,c$ parameterize the family of traveling waves.
Thus the modulational stability index has the following representation
\[
\Delta = \frac{1}{2}\left( \{T,P\}_{E,c} +2 \{M,P\}_{a,E} \right)^3
- 3 \left(\frac{3}{2}\{T,M,P\}_{a,E,c}\right)^2.
\]
\end{thm}

\begin{proof}
Let $w_i(x;\mu)$, $i=1,2,3$, be three linearly independent solutions of \eqref{gkdv}, and let
$\WM(x,\mu)$ be the solution matrix with columns $w_i$.  Expanding the above solutions in powers
of $\mu$ as
\[
w_i(x,\mu)=w_i^0(x)+\mu w_i^1(x)+\mu^2 w_i^2(x)+\mathcal{O}(|\mu|^3)
\]
and substituting them into \eqref{lineargkdv}, the leading order equation
becomes
\[
\frac{d}{dx}w_i^0(x)=\HM(x;0)w_i^0(x).
\]
Using Proposition \ref{Kernel}, we choose $w_i(x)=Y_i(x)$ where the vectors $Y_i(x)$ are defined in equation \eqref{solutionvec}.
The higher order terms in the above expansion yield
\[
\frac{d}{dx}w_i^{j}(x)=\HM(x;0)w_i^{j}(x)+V_i^{j-1}(x),\;\;j\geq 1,
\]
where $V_i^{j-1}=\left(0,0,-(w_i^{j-1})_1\right)^{t}$ and $(v)_1$ denotes the first component of the vector $v$.
Notice that for each of the higher order terms $j\geq 1$, we require $w_j^i(0)=0$.  This implies
that $\WM(0,\mu)=\UM(0,0)$ in a neighborhood of $\mu=0$, where $\UM(0,0)$ is defined in \eqref{MatrixSolutionInitial}.
In the case $j=1$, the $i=1$ equation is equivalent to the equation $L_0 w_1^1=u_x$.  It follows again from
Proposition \ref{Kernel}  that we can choose
\[
w_1^1(x)=\left(
           \begin{array}{c}
             -u_c \\
             -u_{cx} \\
             -u_{cxx} \\
           \end{array}
         \right)
         +
         u_{-}\left(
                \begin{array}{c}
                  u_a \\
                  u_a' \\
                  u_a'' \\
                \end{array}
              \right)
              -\frac{u_{-}^2}{2}
              \left(
                \begin{array}{c}
                  u_E \\
                  u_E' \\
                  u_E'' \\
                \end{array}
              \right).
\]
Notice the above coefficients of $Y_2(x)$ and $Y_3(x)$ are determined by differentiating $E-V(u_{-};a,c)=0$ with respect to the
parameters $a$, $E$, and $c$.  Moreover, using variation of parameters as well as the
identities $\{u_x,u\}_{x,E}=-1$ and $\{u,u_x\}_{x,a}=u$, we choose
\begin{eqnarray}
w_i^j(x)&=&\WM(x,0)\int_0^x \WM(z,0)^{-1}V_i^{j-1}(z)dz\nonumber\\
&=&\left(
     \begin{array}{ccc}
       u_x\int_0^x(w_i^{j-1})_1\{u,u_x\}_{a,E}\;dz-u_a\int_0^x(w_i^{j-1})_1dz+u_E\int_0^x(w_i^{j-1})_1u\;dz \\
       u_{xx}\int_0^x(w_i^{j-1})_1\{u,u_x\}_{a,E}\;dz-u_{ax}\int_0^x(w_i^{j-1})_1dz+u_{Ex}\int_0^x(w_i^{j-1})_1u\;dz\\
       u_{xxx}\int_0^x(w_i^{j-1})_1\{u,u_x\}_{a,E}\;dz-u_{axx}\int_0^x(w_i^{j-1})_1dz+u_{Exx}\int_0^x(w_i^{j-1})_1u\;dz\\
     \end{array}
   \right) \label{variationparameters}
\end{eqnarray}
for $w_i^1$, $i=2,3$, and $w_i^j$ for $i>1$.
Finally, by \eqref{MatrixSolutionDifference}, we have the following expression valid as $\mu\to 0$:
\[
\delta\WM(\mu)=
\left(
  \begin{array}{ccc}
    \mathcal{O}(|\mu|^2) & \mathcal{O}(|\mu|) & \mathcal{O}(|\mu|) \\
    \mu V'(u_{-})\left(-T_c+u_{-}T_a-\frac{u_{-}^2}{2}T_E\right)+\mathcal{O}(|\mu|^2) & V'(u_{-})T_a+\mathcal{O}(|\mu|) & V'(u_{-})T_E+\mathcal{O}(|\mu|) \\
    \mathcal{O}(|\mu|^2) & \mathcal{O}(|\mu|) & \mathcal{O}(|\mu|) \\
  \end{array}
\right),
\]
where $\delta\WM(\mu):=\WM(x,\mu)\big{|}_{x=0}^T$, and the $\mathcal{O}(\mu)$ and $\mathcal{O}(\mu^2)$ terms above are computed using \eqref{variationparameters}.
Note that all of the $\mathcal{O}(\mu,\mu^2)$ terms in the above are
necessary for the calculation, however we do not write them out.
Recalling that our choice of basis implies
$\det(\WM(0,\mu))=-1$,
we have
\begin{eqnarray*}
D(\mu,1)&=&-\textrm{det}\left(\delta\WM(\mu)\right)\\
&=&-\frac{1}{2}\{T,M,P\}_{a,E,c}\;\mu^3+\mathcal{O}(|\mu|^4),
\end{eqnarray*}
from which the expression for $\tr(\MM_{\mu\mu\mu}(0))$ follows by Theorem \ref{realspectrum}.  Moreover,
it follows from Lemma \ref{vanishlem} and the fact $\MM(\mu)=\delta\WM(\mu) \WM(0,0)^{-1}+{\bf I}$ and a rather tedious calculation that
\begin{eqnarray*}
\tr(\MM_{\mu\mu}(0))&=&-2\mu^{-2}\tr(\textrm{cof}(\MM(\mu)-I))\big{|}_{\mu=0}\\
&=&-2\left(-\frac{1}{2}\{T,P\}_{E,c}-\{M,P\}_{E,a}\right)
\end{eqnarray*}
as claimed.
\end{proof}

\begin{corr}
$\{T,M,P\}_{a,E,c}<0$ is a sufficient condition for a non-trivial intersection
of $\spec(\partial_x\mathcal{L}[u])$ with the real axis.
\end{corr}

\begin{proof}
This is now clear from Theorems \ref{realspectrum} and \ref{TraceIdentities}.
\end{proof}

At this point we can make a connection to the stability theory
for the solitary waves.
\begin{corr}
In the case of power-law nonlinearity and wavespeed $c>0$, there are always unstable periodic traveling
waves in a neighborhood of the solitary wave $(a=E=0)$ if $p>4$.  Moreover, such long wavelength
periodic waves exhibit a modulational instability if and only if $p>4$.
\end{corr}
\begin{proof}
First, note that the scaling invariance in equation \eqref{scaling} implies the periodic solution $u(x;a,E,c)$
satisfies
\[
u(x;a,E,c)=c^{1/p}u\left(c^{1/2}x;\frac{a}{c^{1+1/p}},\frac{E}{c^{1+2/p}},1\right),
\]
which allows us to compute $T_c$, $M_c$, $P_c$ explicitly as follows:
\begin{align*}
T_c &= -\frac{1}{2c}T-\frac{a(p+1)}{pc}T_a - \frac{E(p+2)}{pc}T_E \\
M_c &= \left(\frac{1}{pc} - \frac{1}{2c}\right)M + \left(T_c + \frac{T}{2c}\right)u_- -\frac{a(p+1)}{pc}\left(M_a - T_a u_-\right) -\frac{E(2+p)}{pc}\left(M_E - T_E u_-\right)\\
P_c &= \left(\frac{2}{pc} - \frac{1}{2c}\right)P + \left(T_c + \frac{T}{2c}\right)u_-^2 - \frac{a(p+1)}{pc}\left(P_a - T_a u_-^2\right) - \frac{E(2+p)}{pc} \left(P_E - T_E u_-^2\right),
\end{align*}
where $T_c$ follows from equation \eqref{period}.
Since we know that $P_E=2T_c$ and $P_a = 2 M_c$ the above
serves to simplify the last row and column.
When $a$ and $E$ are small there are two turning points $r_1,r_2$ in the
neighborhood of the origin and a third turning point $r_3$ which
is bounded away from the origin. In the solitary
wave limit $a,E\rightarrow 0$ we have $r_1-r_2=O(\sqrt{a^2-2E})$. In this
limit we have the following asymptotics for small $(a,E)$
\begin{align*}
M(a,E,1) &= O(1) \\
P(a,E,1) &= O(1) \\
T(a,E,1) &= O\left(\ln(a^2-2E)\right) \\
T_a(a,E,1) &=O\left(\frac{a}{a^2-2E}\right) = M_E \\
T_E(a,E,1)  &=O\left(\frac{1}{a^2-2E}\right).
\end{align*}
Thus the asymptotically largest minor of $\{T,M,P\}_{a,E,c}$ is $-T_E M_a P_c$, from which it follows
\[
\{T,M,P\}_{a,E,c}\sim -T_E M_a \left(\frac{2}{pc} - \frac{1}{2c}\right) P
\]
as $a,E\to 0$.  It is known that for traveling waves below the separatrix (see appendix)
that under some minor convexity assumptions $T_E>0$. It can also
be shown that (see appendix) $M_a(a,0)<0$ for $E=0$ and $a$ sufficiently
small.  Thus the orientation index $\{T,M,P\}_{a,E,c}$ is negative for
$p>4$ and $a,E$ sufficiently small (in other words sufficiently
close to the solitary wave) and positive for $p<4$ and $a,E$ sufficiently
small.
This also follows, of course, from Gardner's long-wavelength
theory\cite{Gard} but it provides a good check for the present theory.

To prove the second claim, notice the above asymptotics implies
\begin{eqnarray*}
\tr(\MM_{\mu\mu}(0))&\sim& T_E\left(\frac{2}{pc}-\frac{1}{2c}\right)P
\end{eqnarray*}
in the limit as $a,E\to 0$.  Hence, it follows that
$\sgn\Delta(a,E,c)=\sgn(4-p)$ for $a,$ $E$ sufficiently small, and periodic
waves of sufficiently long period are also modulationally unstable for $p>4$.

\end{proof}

\begin{rem}

It is worth noting that the instability mechanism detected
by the discriminant $\Delta$ is not
present in the solitary wave case: in the solitary wave limit the
bands of spectrum connected to the origin collapse to the origin.
This instability does {\em not} appear to follow from Gardner's calculation:
Gardner shows that the point eigenvalue of the solitary wave opens into a
small loop of spectrum, predicting the real eigenvalues
detected by $\{T,M,P\}_{a,E,c},$ but the modulational instability
detected by $\Delta$ is not detected.
Thus suggests the heuristic that periodic solutions should go unstable
before the solitary waves. The small amplitude stability calculation of
H{\v a}r{\v a}gu{\c s} and
Kapitula for the generalized KdV equation amounts to a calculation of
this discriminant in that limiting case, and their proof that the
small amplitude waves go unstable at $p=2$ is the first result we are aware
of along these lines.

We believe that a small amplitude analysis of $\Delta(a,E,c)$ should
be possible.  It follows by a simple calculation that $\Delta=0$
at the stationary solution.  By expanding near by solutions in terms of amplitude
instead of the energy $E$, we believe the first non-zero term of the discriminant
should be proportional to a polynomial which switches signs at $p=2$, thus recovering the
small amplitude result of H\v{a}r\v{a}gu\c{s} and Kapitula \cite{HK}.
We have not as yet carried out such an analysis.
\end{rem}

Using the identities derived in Appendix $1$, we now have a sufficient
criterion for the existence of a non-trivial intersection of
$\spec(\partial_x\mathcal{L}[u])$ with the real axis in terms of the conserved
quantities $M$, $P$ and $H$ of the gKdV flow, as well as a necessary and
sufficient condition for understanding the normal form of the spectrum
in a neighborhood of the origin. It is a rather striking fact that
both of these indices can be expressed entirely in terms of
the conserved quantities of the flow. The monodromy itself depends
on $u_-(a,E,c)$, the classical turning point of the traveling wave, as well
as various functions and derivatives of this quantity, but the
indices themselves only depend on the conserved quantities. This is,
of course, the Whitham philosophy, but we are only aware of a few cases
(other than the integrable calculations, which are very special) in which
make this rigorous.

In the next section we outline the connections of this calculation to
a calculation based more directly on the linearized operator. While
not strictly necessary this calculation is useful since it clarifies the
way in which various bifurcations can occur. In this section we calculate
the null-space and generalized null-space of the linearized operator
and sketch a perturbation calculation analogous to the one given for the
Evan's function.

\section{Local Analysis of $\spec(\partial_x\mathcal{L}[u])$ via the Floquet-Boch Decomposition}

\subsection{Floquet-Bloch Decomposition}

In this section we sketch an approach to this problem working directly
with the linearized operator rather than with the Evan's function.
While these two approaches are presumably equivalent the former
seems less straightforward than the latter. In particular it
is not clear how one might derive the orientation index in
this way, and the calculation of the modulational stability index
gives a quantity which seems much less transparent. Nevertheless we present
an outline of this calculation (omitting some details) since it does give
some insight into the results of the
previous section.

From Floquet theory, we know any bounded eigenfunction $v(x)$
of $\partial_x\mathcal{L}[ u]$ must satisfy
\[
v(x+T)=e^{i\gamma}v(x)
\]
for some $\gamma\in[-\pi,\pi]$.  The quantity $e^{i\gamma}$ is
known as the Floquet multiplier of the eigenfunction $v$.
It follows any eigenfunction $v(x)$ can be represented in the form
$v(x)=e^{i\gamma x/T}P(x)$ where $P(x+T)=P(x)$.
The fact that $\partial_x\mathcal{L}[ u]v(x)=\mu v(x)$ for some $\mu\in\CM$
implies
\[
e^{i\gamma x/T}J_\gamma \mathcal{L}_\gamma [u]P(x)=\mu e^{i\gamma x/T} P(x)
\]
where $J_\gamma=\left(\partial_x+i\frac{\gamma}{T}\right)$ and $\mathcal{L}_\gamma[
u]=-\left(\partial_x+i\frac{\gamma}{T}\right)^2-f'( u)+c$ are the so-called
Bloch operators.  This suggests fixing a $\gamma\in[-\pi,\pi]$ and considering the
eigenvalue problem for the operator $J_\gamma\mathcal{L}_\gamma[u]$
on the Hilbert space $\mathcal{H}=L^2(\RM / T\ZM;\CM)$.  This procedure is
known as a Bloch decomposition of the eigenvalue problem
\eqref{lineargkdv} and we consider the Bloch operators operators as acting
on $\mathcal{H}$.  Notice for $\gamma\neq 0$ the operators $J_\gamma\mathcal{L}_\gamma[u]$ are closed
in this space with compactly embedded domain $H^3(\RM / T\ZM;\CM)$.
It follows that these operators have compact resolvent and hence their
spectra consists of only point
spectra with finite algebraic multiplicities. Moreover, one has
\[
\spec(\partial_x\mathcal{L}[
u])=\bigcup_{\gamma\in[-\pi,\pi]}\spec\left(J_\gamma \mathcal{L}_\gamma[
u]\right).
\]
Thus, this decomposition reduces the problem of locating the
continuous spectrum of the
operator $\partial_x\mathcal{L}[u]$ on $L^2$ to the problem of determining the
discrete spectrum of a one parameter family of operators $\{J_\gamma\mathcal{L}_\gamma[u]\}_{\gamma\in[-\pi,\pi]}$
on $\mathcal{H}$.  Our first goal is to understand the nature of the spectrum of the operator
$J_0\mathcal{L}_0[u]$ at the origin.  Notice in particular that for $\gamma_1,\gamma_2\in[-\pi,\pi]$,
$J_{\gamma_1}\mathcal{L}_{\gamma_1}[u]$ is a compact perturbation of $J_{\gamma_2}\mathcal{L}_{\gamma_2}[u]$, and hence
routine calculations prove the above parameterization of $\spec(J_\gamma\mathcal{L}_\gamma[u])$ is in fact
continuous.  We then consider the operator $J_\gamma\mathcal{L}_\gamma[u]$ for $|\gamma|\ll 1$,
treating it as a small perturbation of $J_0\mathcal{L}_0[u]$, with our end goal
being to study how the spectrum bifurcates from the $\gamma=0$ state.

We begin with analyzing the generalized periodic null space of the operator $\partial_x\mathcal{L}[u]$, denoted
$N_g(\partial_x\mathcal{L}[u])=\bigcup_{n=1}^\infty N((\partial_x\mathcal{L}[u])^n)$.

\begin{prop}\label{gennull}
Suppose that the Jacobians $\{T,M\}_{a,E}$ and $\{T,P\}_{a,E}$ do not
simultaneously vanish, and $\{T,M,P\}_{a,E,c}\neq 0$.
Then zero is an eigenvalue of the operator
$\partial_x\mathcal{L}[u]= J_0\mathcal{L}_0[u]$ considered
on $\mathcal{H}$ of algebraic multiplicity three and geometric multiplicity
two. In the case $\{T,M\}_{a,E}\neq 0$, the functions
\begin{eqnarray*}
\phi_0&=& \{T,u\}_{a,E},\;\;\;\;\;\;\;\;\;\;\;\;\;\;\;\;\;\;\;\;\;\;\;\;\;\psi_0=
1,\\
\phi_1&=& \{T,M\}_{a,E}\;u_x,\;\;\;\;\;\;\;\;\;\;\;\;\;\;\;\;\;\;\psi_1=\int_0^x\phi_2(s)ds,\\
\phi_2&=& \{u,T,M\}_{a,E,c}\;\;\;\;\;\;\;\;\;\;\;\;\;\;\;\;\;\;\;\psi_2=\{T,M\}_{E,c}+\{T,M\}_{a,E}u,
\end{eqnarray*}
provide a basis for $N_g(\partial_x\mathcal{L}[u])$ and
$N_g(\mathcal{L}[u]\partial_x)$ respectively.
Specifically we have the relations
\begin{eqnarray*}
\partial_x\mathcal{L}[u] \phi_0&=&0~~~~~~~~~~~~~~~~~~~~~~~~~~~\mathcal{L}[u]\partial_x \psi_0 = 0 \\
\partial_x\mathcal{L}[u] \phi_1&=&0~~~~~~~~~~~~~~~~~~~~~~~~~~~\mathcal{L}[u]\partial_x \psi_1 = -\psi_2 \\
\partial_x\mathcal{L}[u] \phi_2&=&-\phi_1~~~~~~~~~~~~~~~~~~~~~~~\mathcal{L}[u]\partial_x \psi_2 = 0.
\end{eqnarray*}
\end{prop}

\begin{proof}
The constants above are chosen for convenience, and the functions
above are not normalized. For instance, $\phi_2$ can be any multiple
of $u_x$ and similarly $\psi_0$ any constant.  Also, the ordering
is chosen so that $\left<\phi_j,\psi_k\right>=0$ for $i\neq k$.
Notice this proposition does not follow directly from Proposition
\ref{Kernel} since the functions $ u_a$, $ u_E$ and $ u_c$ are
not in general T-periodic, and one must chose linear combinations
which are periodic and thus belong to $\mathcal{H}.$

First observe that \eqref{MatrixSolutionDifference} implies
$\phi_0$ and $\phi_1$ are $T$-periodic and
belong to $N(J_\gamma \mathcal{L}_\gamma[u])$. In particular, Corollary
\ref{NullSpaces} from the appendix implies
$N(\partial_x\mathcal{L}[u])=\textrm{span}\{\phi_0,\phi_1\}$ when
considered as an operator on $\mathcal{H}$. The fact that the
monodromy at the origin is the identity plus a rank one perturbation
suggests that there are two linear combinations that can be chosen
to be periodic.  Specifically
we define
\[
\phi_2=\left|%
\begin{array}{ccc}
  u_a & T_a & \int_0^T u_a dx \\
  u_E & T_E & \int_0^T u_E dx \\
  u_c & T_c & \int_0^T u_c dx \\
\end{array}%
\right| = \{u,T,M\}_{a,E,c}
\]
and it is clear from \eqref{MatrixSolutionDifference} that
$\phi_2\in \mathcal{H}$ and $J_0\mathcal{L}_0[u]\phi_2=\phi_1$ as claimed.
Thus, if $\{T,M\}_{a,E}\neq 0$, $\phi_2$ gives
a function in $N((J_0\mathcal{L}_0[u])^2)- N(J_0\mathcal{L}_0[u])$.

Similarly, Corollary \ref{NullSpaces} implies that $\psi_0$ and $\psi_2$
are belong to $N(\mathcal{L}[u]\partial_x)$, and are linearly independent
provided that $\{T,M\}_{a,E}\neq 0$.  Moreover, a it is clear from construction
that $\psi_1\in\mathcal{H}$ and a straightforward computation shows that
$\psi_1$ belongs to $N((\mathcal{L}_0[u]J_0)^2)-N(\mathcal{L}_0[u]J_0)$ as claimed.

In order to complete the proof, we must now show these three functions comprise
the entire generalized null space of $J_0\mathcal{L}_0[u]$ on $\mathcal{H}$.  To this end, we prove that
neither of the functions $\phi_0,\phi_2$ belong to the range of $J_0\mathcal{L}_0[u]$
by appealing to the Fredholm alternative.  It follows that the equation $J_0\mathcal{L}_0[u] v=\phi_0$ has a solution
in $\mathcal{H}$ if and only if the following solvability conditions are
simultaneously satisfied:
\begin{eqnarray*}
\left<1,\phi_0\right>&=&\{T,M\}_{a,E}=0,\;\;\textrm{and}\\
\left<u,\phi_0\right>&=&\frac{1}{2}\{T,P\}_{a,E}=0.
\end{eqnarray*}
Thus, if either $\{T,M\}_{a,E}$ or $\{T,P\}_{a,E}$ are non-zero, then
$N((J_0\mathcal{L}_0[u])^2)-N(J_0\mathcal{L}_0[u])=\textrm{span}\{\phi_2\}$.  Similarly,
$N((J_0\mathcal{L}_0[u])^3)-N((J_0\mathcal{L}_0[u])^2)\neq\emptyset$ if and only if
the equation $L_0 v=\phi_2$ has a solution in $\mathcal{H}$, i.e. if and only if
\begin{eqnarray*}
\left<u,\phi_2\right>=\frac{1}{2}\{T,M,P\}_{a,E,c}=0,
\end{eqnarray*}
which finishes the proof.

A similar construction in the case $\{T,M\}_{a,E}=0$ but $\{T,P\}_{a,E}\neq 0$
gives a basis in this case.
\end{proof}

\begin{rem}
It is worth remarking in some detail on the physical significance of these
conditions and the relationship to the Whitham modulation theory. Obviously
$(a,E,c)$ are constants of integration arising in the ordinary differential
equation defining the traveling wave, and $T,M,P$ are constants of the
PDE evolution. One of the main ideas of the Whitham theory is to locally parameterize the wave by the constants of motion. The non-vanishing of the
Jacobians is exactly what allows one to do this. Non-vanishing of
$\{T,M,P\}_{a,E,c}$ is  equivalent to demanding that locally
the map
$(a,E,c) \mapsto (T,M,P)$ have a unique $C^1$ inverse - in other words the
conserved quantities $(T,M,P)$ are good local coordinates for the family
of traveling waves. Similarly non-vanishing of one of $\{T,M\}_{a,E}$ and $\{T,P\}_{a,E}$ is, at least for periodic waves below the separatrix, equivalent to
demanding that the matrix
\[
\left(\begin{array}{ccc} T_a & M_a & P_a \\ T_E & M_E & P_E \end{array}\right)
\]
have full rank, which is  equivalent to demanding that the map
$(a,E) \mapsto (T,M,P)$ (at fixed $c$) have a unique $C^1$ inverse - in other
words two of the conserved quantities
give a smooth parameterization of the family of  traveling waves of
fixed wave-speed.
As long as $E\neq 0$ we can use the identities developed in the appendix
to eliminate $T$ in favor of $H$. Thus in the case $E\neq 0$ (which {\em does not} include the solitary wave wave) the null-space being two dimensional is
equivalent to  two of the conserved quantities $(M,P,H)$ giving
a $C^1$ parameterization
of the traveling wave solutions at constant wavespeed, and the space
$N((J_0\mathcal{L}_0[u])^2)-N(J_0\mathcal{L}_0[u])$ being one dimensional
is equivalent to the three conserved quantities $(M,P,H)$ giving
 a $C^1$ parameterization of the full family of traveling waves.
\end{rem}

Notice it follows the vanishing of $\{T,M,P\}_{a,E,c}$ is connected with a change in the Jordan structure
of the linearized operator $J_0\mathcal{L}_0[u]$ considered on $\mathcal{H}$: $\{T,M,P\}_{a,E,c}\neq 0$
ensures the existence of a non-trivial Jordan piece in the generalized null space of 
of dimension exactly one.
Moreover, it guarantees that the variations in the constants associated to the
family of traveling wave solutions by reducing \eqref{gkdv} to quadrature are
enough to generate the entire generalized periodic null space of the operator $J_0\mathcal{L}_0[u]$.
Henceforth, we shall assume $\{T,M,P\}_{a,E,c}\neq 0$ and that
$\{T,M\}_{a,E}\neq 0$ - trivial modifications are necessary
if $\{T,M\}_{a,E}$ vanishes but $\{T,P\}_{a,E}$ does not.

\subsection{Analyticity of Eigenvalues Bifurcating from $\mu=0$}

Our next goal then is to consider the operator $J_\gamma\mathcal{L}_\gamma[u]$ for small $\gamma$,
treating it as a small perturbation of $J_0\mathcal{L}_0[u]$.  To this end, notice that if we define
$L_0:=J_0\mathcal{L}_0[u]$, $L_1:=\mathcal{L}_0[u]-2\partial_x^2$, and $L_2:=-3\partial_x$, it follows
that
\[
J_\gamma\mathcal{L}_\gamma[u]=L_0+\varepsilon L_1+\varepsilon^2L_2-\varepsilon^3,
\]
where $\varepsilon$ is related to the Floquet exponent via $\varepsilon=\frac{i\gamma}{T}$.
By Proposition \ref{gennull}, we know the operator $L_0$ has three periodic
eigenvalues at the origin.  Our present goal is to determine how these eigenvalues
bifurcate from the $\gamma=0$ state.  In this section we only sketch the
relevant details - for similar calculations see the papers of
Ivey and Lafortune\cite{IL}, or Kapitula, Kutz and Sanstede.\cite{KKS}

Since the Hilbert space $\mathcal{H}$
consists of T-periodic functions, eigenvalues of $J_\gamma\mathcal{L}_\gamma[u]$ correspond
to $1$ being an eigenvalue of the monodromy operator $\Phi(T;\mu,\varepsilon)$ for to the eigenvalue
problem $J_\gamma\mathcal{L}_\gamma[u]v=\mu v$.  Thus, it is natural to introduce the following
``modified'' periodic Evans function
\[
D_0(\mu,\varepsilon) = \textrm{det}\left(\Phi(T;\mu,\varepsilon)-I\right).
\]
Notice that $D_0(\mu,\varepsilon)$ is clearly an analytic function of the two complex variables $\mu$ and $\varepsilon$.
Our first goal then is to analyze the possible behavior of of the solutions of $D_0(\mu,\varepsilon)=0$
in a small neighborhood of $(0,0)$.

\begin{lem}\label{evprep}
Let $F(x,y)$ be a complex valued function of two complex variables $x$ and $y$ which is analytic
in a neighborhood of $(0,0)\in\CM^2$.  Moreover, suppose that $F(0,0)=F_x(0,0)=F_{xx}(0,0)=0$, $F_{xxx}(0)\neq 0$,
and $F_y(0,0)=0$.  Then for small $y$, the equation $F(x,y)=0$
has three roots in a neighborhood of the origin.  Moreover, these roots are given by $(x,y)=(f_j(y),y)$, $j=1,2,3$,
where the $f_j$ satisfy one of the following conditions:
\begin{enumerate}
  \item[(i)] One function $f_j$ can be expressed as a Puiseux series as $f_j(y)=\sum_{n=1}^\infty a_n^j y^{n/2}$
              in a neighborhood of $y=0$, where $a_1\neq 0$.
  \item[(ii)] Two of the functions $f_j$ admits a Puiseux series representation of the form $f_j(y)=\sum_{n=2}^\infty a_n^j y^{n/3}$
              in a neighborhood of $y=0$, where $a_2\neq 0$.
  \item[(iii)] All three functions $f_j$ are $\mathcal{O}(\varepsilon)$ and are analytic in $y$ in a neighborhood of $y=0$, i.e.
           they can be represented as $f_j(y)=\sum_{n=1}^\infty a_n^j y^n$ where $a_1\neq 0$ assuming $F_{yyy}(0,0)\neq 0$.
\end{enumerate}
In the case \textit{(iii)}, if $F_{yyy}(0,0)=0$ then all three eigenvalues are analytic in $\varepsilon$, with two eigenvalues
of order $\mathcal{O}(|\varepsilon|)$ and the remaining eigenvalue of order at least $\mathcal{O}(|\varepsilon|^2)$.
\end{lem}

\begin{proof}
By the Weierstrass preparation theorem, the function $F(x,y)$ can be expressed as
\[
F(x,y)=\left(x^3+\eta_2(y)x^2+\eta_1(y)x+\eta_0(y)\right)h(x,y)
\]
for small $x$ and $y$, where each $\eta_j$ is analytic, and $h$ is analytic satisfying $h(0,0)\neq 0$.
It follows the three roots of $F(x,y)$ near $(0,0)$ are determine by the cubic polynomial
$G(x,y)=x^3+\eta_2(y)x^2+\eta_1(y)x+\eta_0(y)$.  By hypothesis, we have that
$\eta_j(0)=0$ for $j=0,1,2$, $\eta_0'(0)=0$, and $\eta_0'''(0)\neq 0$.  Hence, the Newton diagram
for the equation $G(x,y)=0$ is the same as that in figure 1, from which the lemma follows.
\end{proof}

We now wish to apply Lemma \ref{evprep} to the equation $D_0(\mu,\varepsilon)=0$, with $x=\mu$ and $y=\varepsilon$,
and use the Fredholm alternative to show only possibility $\textit{(iii)}$ can occur.  Notice that Theorem \ref{TraceIdentities} implies
$\frac{\partial^k}{\partial \mu^k} D_0(\mu;0)=0$ for $k=0,1,2$ and, moreover, $\frac{\partial^3}{\partial\mu^3}D_0(\mu;0)\neq 0$
under the assumption $\{T,M,P\}_{a,E,c}\neq 0$.  To apply Lemma \ref{evprep} then, we need the following lemma.

\begin{lem}
We have $\frac{\partial}{\partial\varepsilon}D_0(0,0)=0$.
\end{lem}

\begin{proof}
This proof proceeds much like that of Theorem \ref{TraceIdentities}.  Defining $W(x;\mu,\varepsilon)$
to be the solution matrix to the first order system corresponding to $J_\gamma\mathcal{L}_\gamma[u]v=\mu v$
written in the basis $Y_i(x)$ defined in \eqref{solutionvec}, arguments similar to those above yield
for small $\varepsilon$
\[
\textrm{det}\left(W(T;0,\varepsilon)-W(0;0,\varepsilon)\right)=
\left(
  \begin{array}{ccc}
    \mathcal{O}(|\varepsilon|) & \mathcal{O}(|\varepsilon|) &              \mathcal{O}(|\varepsilon|) \\
    \mathcal{O}(|\varepsilon|) & V'(u_{-})T_a+\mathcal{O}(|\varepsilon|) & V'(u_{-})T_E+\mathcal{O}(|\varepsilon|) \\
    \mathcal{O}(|\varepsilon|) & \mathcal{O}(|\varepsilon|) &           \mathcal{O}(|\varepsilon|)\\
  \end{array}
\right),
\]
and hence $D_0(0,\varepsilon)=\mathcal{O}(|\varepsilon|^2)$ as claimed.
\end{proof}

We are now in position to prove our main result of this section.  By the
above work, we can apply Lemma \ref{evprep} to the equation
$D_0(\mu,\varepsilon)=0$.  The next theorem uses the Fredholm alternative to
discount possibilities \textit{(i)} and \textit{(ii)} from emma
\ref{evprep}, and establish the analyticity of the eigenvalues near $\mu=0$.
Basically this amounts to checking that (generically) the null-space of the
linearized operator has the same Jordan structure as the monodromy
map at the origin.

\begin{thm}\label{analyticev}
For small $\varepsilon$, the linear operator $J_\gamma\mathcal{L}_\gamma[u]$
has three eigenvalues which bifurcate from $\mu=0$ and are analytic in $\varepsilon$.
\end{thm}

\begin{proof}
The idea of the proof is to systematically discount possibilities \textit{(i)} and
\textit{(ii)} from Lemma \ref{evprep}, thus leaving only the third possibility.
First, suppose case \textit{(i)} holds.  It follows from the Dunford
calculus that we can expand the eigenvalues and eigenfunctions as
\[\left\{%
\begin{array}{ll}
    \mu=\varepsilon^{1/2}\nu_1+\varepsilon\nu_2+\mathcal{O}(|\varepsilon|^{3/2})\\
    v=f_0+\varepsilon^{1/2}f_1+\varepsilon f_2+\mathcal{O}(|\varepsilon|^{3/2})
\end{array}%
\right.
\]
We will show the assumption that $\{T,M,P\}_{a,E,c}\neq 0$ implies $\nu_1=0$, which yields the desired
contradiction. Using the above expansions of $v$, $\mu$ and $J_\gamma \mathcal{L}_\gamma[\phi]$ in terms of
$\varepsilon$, the leading order equation becomes $L_0f_0=0$.
Thus, $f_0=b_0\phi_0+b_1\phi_1$ for some $b_0,\;b_1\in\CM$.
Continuing, the $\mathcal{O}(|\varepsilon|^{1/2})$ equation turns out to be
$L_0 f_1=\nu_1f_0$. Suppose $\nu_1\neq 0$. By the Fredholm
alternative, this equation is solvable in $\mathcal{H}$ if and only if
$b_0\phi_0+b_1\phi_1\perp N(L_0^\dag)$.  Clearly, $\phi_1\perp
N(L_0^\dag)$ since $\phi_1\in Range(L_0)$.  Moreover, by Lemma \ref{gennull}
$\phi_0\notin N(L_0^\dag)^\perp$ and hence we must have
$b_0=0$ and, with out loss of generality, we take $b_1=1$.
It follows that $f_1$ must satisfy the equation
\[
L_0f_1=\nu_1 \phi_1,
\]
i.e. $f_1=\nu_1 \phi_2+b_2\phi_0+b_3\phi_1$ for some
constants $b_2,\;b_3\in\CM$.

Continuing in this fashion, the $\mathcal{O}(|\varepsilon|)$ equation
becomes
\[
L_0f_2=\nu_1 f_1+\nu_2 f_0-L_1 f_0.
\]
By the Fredholm alternative, this is solvable if and only if
\begin{eqnarray*}
\left<\psi_0,\nu_1 f_1+\nu_2 f_0-L_1
f_0\right>&=&0\;\;\textrm{and}\\
\left<\psi_2,\nu_1 f_1+\nu_2 f_0-L_1 f_0\right>&=&0.
\end{eqnarray*}
By above, $f_0$ is an odd function and since $L_1$ preserves
parity, the solvability condition implies we must require
$\left<\psi_0,f_1\right>=\left<\psi_2,f_1\right>=0$.  However,
this is a contradiction since
$\left<\psi_2,\phi_2\right>=\frac{1}{2}\{T,M\}_{a,E}\{T,M,P\}_{a,E,c}$ and hence it must be that
$\nu_1=0$ as claimed.  Thus, possibility \textit{(i)} can not occur.

Next, assume case \textit{(ii)} of Lemma \ref{evprep} holds.  Then
the Dunford calculus again implies the eigenvalues and eigenvectors can
be expanded in a Puiseux series of the form
\[
\left\{
  \begin{array}{ll}
    \mu=\omega_1\varepsilon^{2/3}+\omega_2\varepsilon^{4/3}+\mathcal{O}(|\varepsilon|^{2}), \\
    v=w_0+\varepsilon^{2/3}w_1+\varepsilon^{4/3}w_2+\mathcal{O}(|\varepsilon|^2).
  \end{array}
\right.
\]
Our goal again is to prove the assumptions $\{T,M,P\}_{a,E,c}\neq 0$ and $\{T,M\}_{a,E}\neq 0$
imply $\omega_1=0$.  Substituting these expansions into $J_\gamma\mathcal{L}_\gamma[u]v=\mu v$ as before,
the leading order equation leads to $w_0=a_0\phi_0+a_1\phi_1$ and the $\mathcal{O}(|\varepsilon|^{2/3})$
equation implies $a_0=0$.  Without loss of generality, we assume $a_1=1$, so that it follows
that $w_1=\omega_1\phi_2+a_2\phi_0+a_3\phi_1$.  The solvability condition at
$\mathcal{O}(|\varepsilon|^{4/3})$ implies that
\[
-\omega_1^2\left<\psi_2,\phi_2\right>=0,
\]
which implies $\omega_1=0$ as above.  Thus, case \textit{(ii)} of Lemma \ref{evprep} can not occur
leaving only case \textit{(iii)}, which completes the proof.
\end{proof}

\subsection{Perturbation Analysis of $\spec(J_\gamma\mathcal{L}_\gamma[u])$ near $(\mu,\gamma)=(0,0)$}

We are now set to derive a modulational stability index in terms of the conserved quantities
of the gKdV flow.  By Theorem \ref{analyticev}, it follows that the eigenvalues
and eigenvectors are analytic in $\varepsilon$, and hence admit a representation of the form
\[
\left\{
  \begin{array}{ll}
    v=v_0+v_1\varepsilon+v_2\varepsilon^2+\mathcal{O}(|\varepsilon|^3), \\
    \mu=\lambda_1\varepsilon+\lambda_2\varepsilon^2+\mathcal{O}(|\varepsilon|^3).
  \end{array}
\right.
\]
At this point, it is tempting to use the functionals
$P_j:=\left<\psi_j,\cdot\right>$ to compute the matrix
action of the operator $J_\gamma\mathcal{L}_\gamma[u]$ onto the corresponding
spectral subspace associated with $N_g(L_0)$.  This would convert
the above eigenvalue problem for a fixed $\gamma$ to the problem
of solving the polynomial equation
\[
\det\left[M_0+\varepsilon M_1+\varepsilon^2 M_2+\mathcal{O}(\varepsilon^3)-\lambda P\right]=0,
\]
at $\mathcal{O}(\varepsilon^2)$, where $M_k=\{P_iL_k\phi_j\}_{i,j}$ and $P=\{P_i\phi_j\}_{i,j}$.
Although this approach has been used to determine stability in the case
where the underlying periodic waves are small (see \cite{GH2} and \cite{HK}),
this approach is flawed in  the current case since, as shown below, the
eigenvector $v$ has a non-trivial projection onto $N_g(L_0)^\perp$
of size $\mathcal{O}(\varepsilon)$.  Since we have no information
about what such a projection would look like, it is unlikely that
one can determine the nature of the spectrum near $\mu=0$ by
computing the matrix action of the operator
$J_\gamma\mathcal{L}_\gamma[u]$ on $\mathcal{H}$ for a general periodic
solution of \eqref{gkdv}.  Instead, we proceed below by developing
a perturbation theory for such a degenerate eigenvalue problem based
on the Fredholm alternative.

Substituting the analytic representation of the eigenvector and
eigenvalue into the equation $J_\gamma\mathcal{L}_\gamma[u] v=\mu
v$, the leading order equation implies $v_0\in N(L_0)$, i.e.
$v_0=c_0\phi_0+c_1\phi_1$ for some $c_0,\;c_1\in\CM$.  At
$\mathcal{O}(|\varepsilon|)$, we get the equation $L_0v_1=(\lambda_1-L_1)v_0$, which
has corresponding solvability conditions
\begin{eqnarray*}
    0=\left<\psi_0,L_0v_1\right>&=&\lambda_1c_0\left<\psi_0,\phi_0\right>-c_0\left<\psi_0,L_1\phi_0\right>
                       - c_1\left<\psi_0,L_1\phi_1\right>,\;\textrm{and} \\
    0=\left<\psi_2,L_0v_1\right>&=&-c_0\left<\psi_2,L_1\phi_0\right> - c_1\left<\psi_2,L_1\phi_1\right>.
\end{eqnarray*}
It follows that we must require $c_0=0$.  Indeed, from the parity
relation $\left<\psi_i,L_k\phi_j\right>=0$ if $i+j+k=0$ mod(2), and the relations
$\left<\psi_0,\phi_0\right>=\{T,M\}_{a,E}\neq
0$ and $\left<\psi_0,L_1\phi_0\right>=T_E\neq 0$, we either have
$c_0=0$ or all three eigenvalues bifurcating from $\mu=0$
have the same leading order  non-zero real part, which is not allowed
by the Hamiltonian symmetries of the spectrum (recall $\spec(\partial_x\mathcal{L}[u])$ is
symmetric about the real and imaginary axis).
With out loss of generality, we then set $c_1=1$ and fix the normalization
\[
\left<\psi_1,v\right>=\left<\psi_1,v_0\right>=-\frac{1}{2}\{T,M\}_{a,E}\{T,M,P\}_{a,E,c}
\]
for all $\varepsilon$.  It follows that $v_0=\phi_1$ and hence $v_1$ satisfies the equation
\[
L_0v_1=(\lambda_1-L_1)\phi_1.
\]
Notice that $L_1v_0=-2\{T,M\}_{a,E}u_{xxx}$
does not belong to $N_g(L_0)$, and hence the eigenfunction $v$ has
a non-trivial projection onto $N_g(L_0)^\perp$ of size
$\mathcal{O}(\varepsilon)$, as claimed above.

We now define $L_0^{-1}$ on $R(L_0)$ with the requirement
that $R(L_0^{-1})$ is orthogonal to $span\{\psi_0,\psi_1\}$.  This requirement ensures that $L_0^{-1}f$ is well-defind
and unique for all $f\in R(L_0)$.  In particular,  it allows us to compute
the projection of $L_0^{-1}f$ onto $N(L_0)$ for each $f\in R(L_0)$.
In order to express the explicit dependence of $v_1$ on $\lambda_1$, we now write
\begin{equation}
v_1=L_0^{-1}\left(\lambda_1-L_1\right)\phi_1+c_2\phi_0+c_3\phi_1\label{v1}
\end{equation}
for some $c_1$, $c_3\in\CM$.  The above normalization condition implies $\left<\psi_1,v_1\right>=0$, i.e.
\[
0=\left<\psi_1,L_0^{-1}(\lambda_1-L_1)\phi_1\right>+c_3\left<\psi_1,\phi_1\right>.
\]
It follows $c_3=0$ by the definition of $L_0^{-1}$ and the fact that
$\left<\psi_1,\phi_1\right>\neq 0$.

Continuing, the $\mathcal{O}(|\varepsilon|^2)$ equation is
\[
L_0v_2=-L_1v_1-L_2v_0+\lambda_1v_1+\lambda_2v_0
\]
with corresponding solvability conditions
\begin{eqnarray*}
0&=&-\left<\psi_0,L_1v_1\right>-\left<\psi_0,L_2v_0\right>+\lambda_1\left<\psi_0,v_1\right>,\;\textrm{and} \\
0&=&-\left<\psi_2,L_1v_1\right>-\left<\psi_2,L_2v_0\right>+\lambda_1\left<\psi_2,v_1\right>.
\end{eqnarray*}
Using the explicit dependence of $v_1$ on $\lambda_1$ and $c_2$, it follows that we can
express the above solvability conditions as
\begin{eqnarray*}
    P_1(\lambda_1)+\widetilde{P}_1(\lambda_1)c_2&=&0\;\textrm{and} \\
    P_2(\lambda_1)-P_0(\lambda_1)c_2&=&0.
\end{eqnarray*}
As this is an over determined system of linear equations for $c_2$, the consistency condition
\[
P(\lambda_1):=P_0(\lambda_1)P_1(\lambda_1)+\widetilde{P}_1(\lambda_1)P_2(\lambda_1)=0
\]
must hold.  In particular, this expresses
$\lambda_1$ as a root of a cubic polynomial with real coefficients.  Since $\varepsilon$
is purely imaginary, modulational stability follows if and only if $P(\lambda)$ has
three real roots, and hence it must be that $\Delta(f;u)$ is a positive multiple of the discriminant of
the cubic polynomial $P(\lambda)$.  Notice that one can explicitly calculate $P(\lambda)$ for a general
non-linearity using just the definitions of the $\phi_j$ and $\psi_j$, except for the inner products
$\left<\psi_0,L_1L_0^{-1}L_1\phi_1\right>$ and $\left<\psi_2,L_1L_0^{-1}L_1\phi_1\right>$.  The first
of these can be calculated regardless of the nonlinearity, but we must restrict to power-law
nonlinearity for the computation of the second(see appendix).
It follows that we can explicitly write down the compatibility condition $P(\lambda_1)=0$
only in terms of the underlying periodic solution $u$ and terms built up out of the generalized
null spaces of $L_0$ and $L_0^{\dag}$ acting on $L^2(\RM/T\ZM)$.  Since the roots of this polynomial
determine the structure of $\spec(J_\gamma\mathcal{L}_\gamma[u])$ in a neighborhood of the origin, we have
proven the following theorem.

\begin{thm}
The periodic solution $u=u(x;a,E,c)$ of \eqref{gkdv} is spectrally unstable in a neighborhood of the origin
if and only if the discriminant $\Delta(a,E,c)$ of the real cubic polynomial
$P(\lambda)$ is positive. Recall that the discriminant of a cubic
$P(\lambda) = a \lambda^3 + b \lambda^2 + c \lambda + d$ is given by
$\Delta = b^2 c^2 - 4 a c^3 - 4 b^3 d - 27 a^2 d^2 + 18 a b c d$
\end{thm}

\begin{rem}
The above result gives a second characterization of
 the modulational stability of
periodic solutions to the generalized Korteveg-DeVries equation with power
law nonlinearity since it is expressed entirely in terms of $T,M,P,H$ and their derivatives, which in turn can be written as functions of $a,E,c$ via integral type formulae. (These are hyperelliptic integrals in the case that $p$ is rational). The formulae remain, however, somewhat daunting. Since this detects the
same instability that the Evan's function based criterion does this quantity
must have the same sign as the discriminant derived in that section, although
we have not shown this explicitly.
\end{rem}

\section{Concluding Remarks}
\subsection{Discussion}

We'd like to consider a concrete example to illustrate our results.
We have chosen to consider the power law gKdV with $p=5 (and c=1)$. In this
case the solitary wave is unstable and hence (by Gardner's result,
which we have checked in this case using our methods) periodic
waves of sufficiently long period are also unstable.
H{\v a}r{\v a}gus and Kapitula\cite{HK} have done some very nice experiments
on this
case using the SpectrUW\cite{SpecTrue,SpecTrue2} package, which they have been
kind enough to share with us. For clarity we have illustrations
representing the spectra they computed numerically, rather than
reproducing their figures - see Figure 3.

\begin{figure}
\begin{center}
\begin{tabular}{ccc}
\includegraphics[scale=.475]{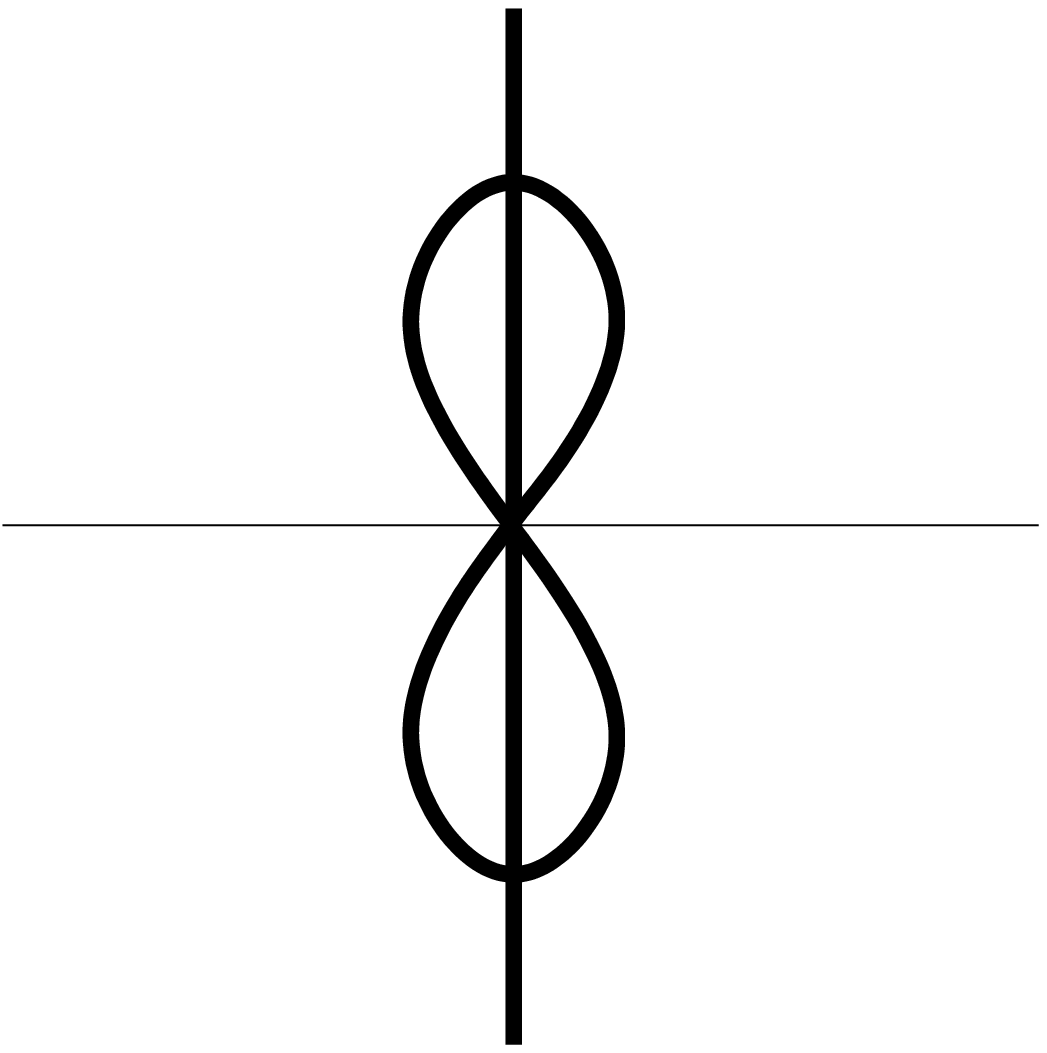}&
\includegraphics[scale=.475]{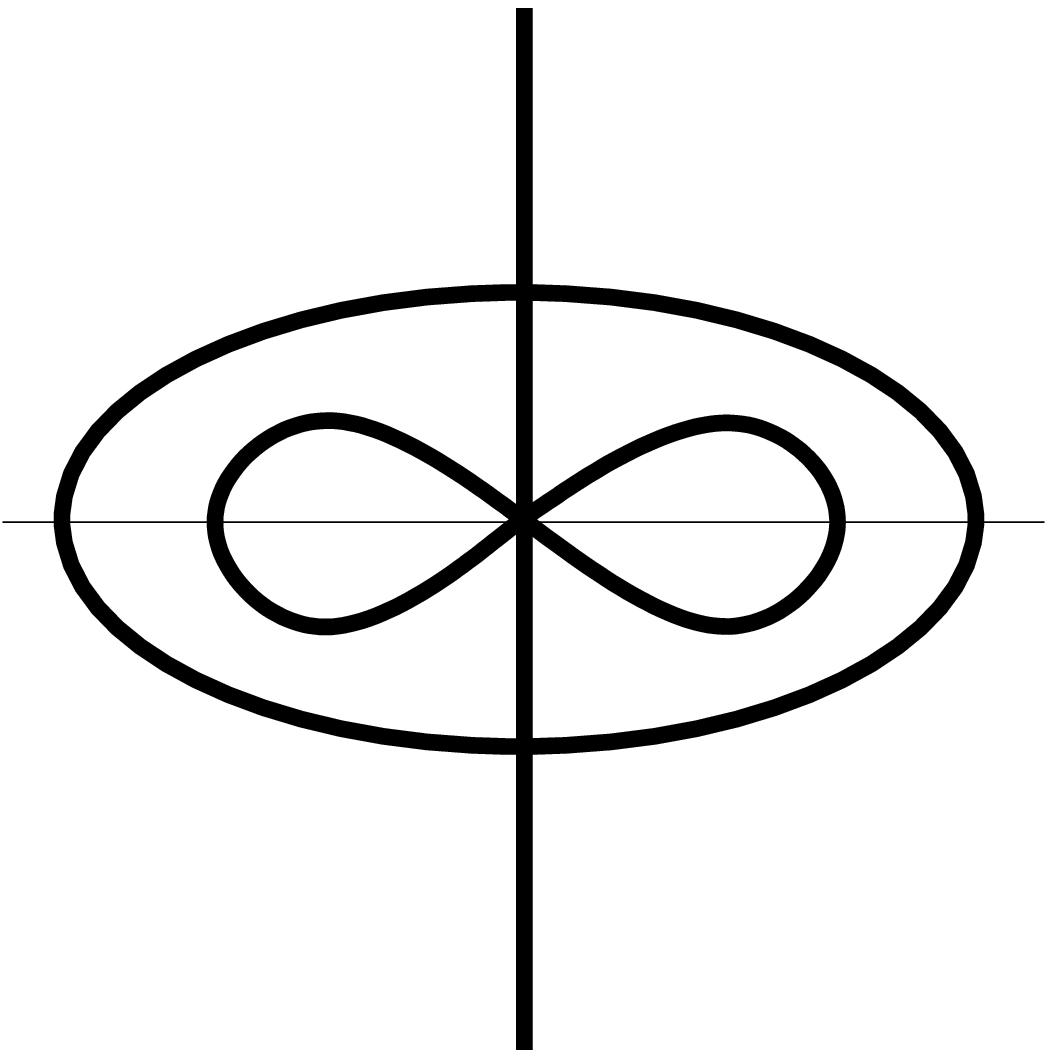}&
\includegraphics[scale=.475]{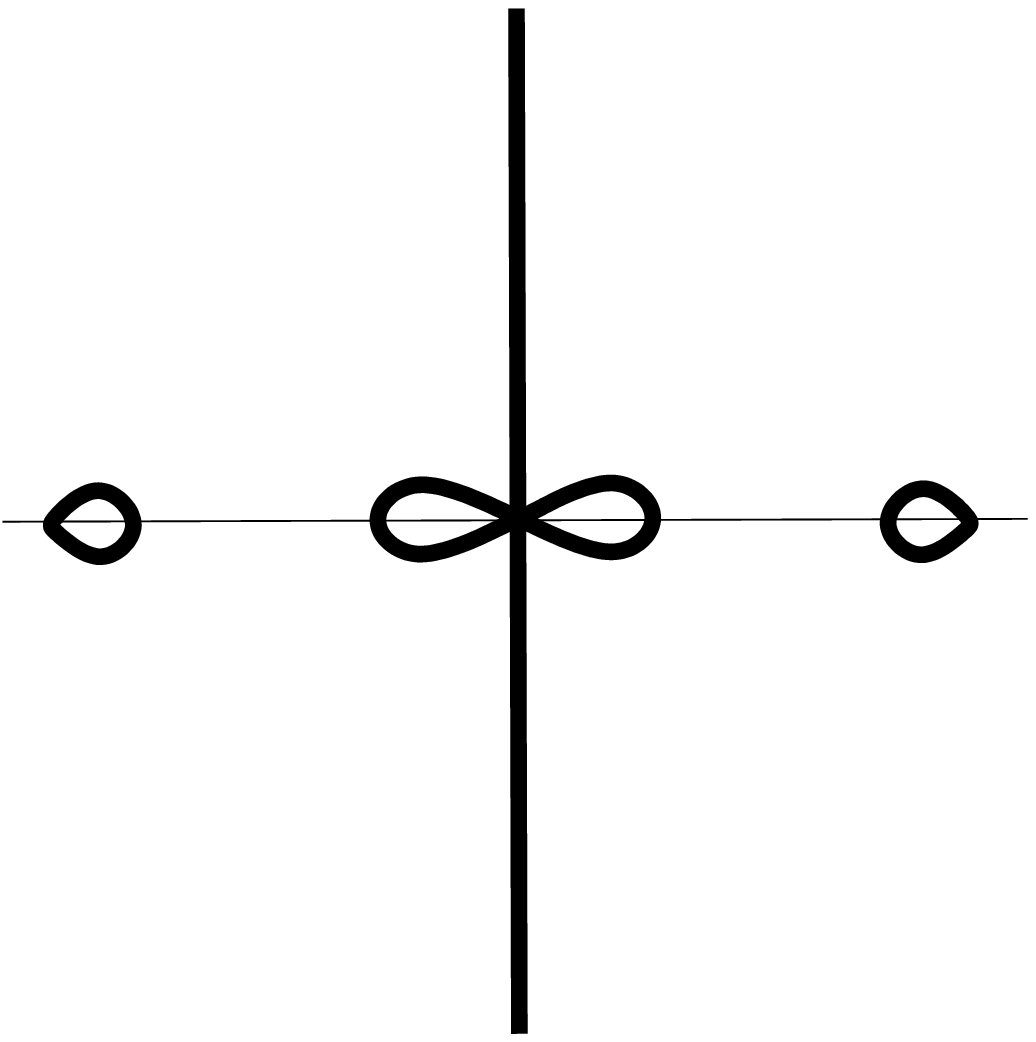} \\
\end{tabular}
\end{center}
\caption{ Cartoon of the spectrum of the linearization of gKdV about a periodic
traveling wave for $p=5,a=0,E=0$ and three different values of the period
(ordered by increasing period). }

\label{fig:goodkrein}
\end{figure}

The first graph in Figure 3 depicts the spectrum for small amplitude periodic
waves (this is the solution branch inside the separatrix). The modulational
instability index $\Delta <0$ indicating a modulational instability,
while the orientation index $\{T,M,P\}_{a,E,c}>0.$ The latter
indicates that the number of eigenvalues on the real axis away from
the origin is even. In this case there are none. The spectrum near the
axis looks like a union of three straight lines, as predicted by the
fact that the normal form of the periodic Evan's function is a homogeneous
polynomial of degree three.  Globally the
spectrum looks like the union of the imaginary axis with a figure eight
shaped curve.

As the period increases one sees spectra which resemble the
second figure, where there is a modulational instability
together with a pair of eigenvalues along the real axis.
 In this case we are still in the case $\Delta <0$,
indicating a modulational instability, and  $\{T,M,P\}_{a,E,c}>0$ indicating
an even number of eigenvalues along the positive real axis. The fact that
these two very different spectral pictures have the same orientation and
modulational instability indices shows that these quantities are
not enough to say qualitatively what the spectral picture looks like,
even in this very simple problem with only one free parameter (the period).

As the period increases still further one sees spectral pictures which
resemble the third picture. As in the previous figure there is an $\infty$
shaped curve of spectrum connected to the origin indicating a modulational
instability ($\Delta<0$) as well as two loops of spectrum intersecting the
real axis and supported away from the origin. These loops are those predicted
by Gardner in his paper arising from the discrete eigenvalues of the
solitary wave problem.
As the period increases and the periodic solution approaches the solitary
wave the circle collapses to a point and the $\infty$ curve collapses to
the origin. The size of both of these features is exponentially small
in the period. In the paper of Kapitula and Haragus the $\infty$ curve
is not visible at the scale of the graph, but it is visible in
numerics they performed for smaller values of period.

Since there is an odd number of eigenvalues on the
real axis in this case (one periodic, two antiperiodic) the
orientation index must now be negative: $\{T,M,P\}_{a,E,c}<0.$ The general
mechanism by which this must occur is clear: a periodic eigenvalue
moves down the real axis, collides with the origin (changing the
Jordan structure of the null-space of the linearized operator, which is
again signalled by the vanishing of  $\{T,M,P\}_{a,E,c}<$) and moves off along
the real axis. However the exact way in which this occurs is not quite clear.

\subsection{Open Problems and Concluding Remarks}
We have considered the stability of periodic traveling waves solutions to
the generalized Korteveg-DeVries equation to perturbations of arbitrary
wavelength. We introduce two indices related to the stability of
the period waves. The first, which is given by the Jacobian of the map
between the constants of integration of the traveling wave ordinary
differential equation and the conserved quantities of the partial differential
equation, serves to count (modulo $2$) the number of periodic
eigenvalues along the real axis. This is, in some sense, a natural
generalization of the analogous calculation for the solitary wave solutions,
and reduces to this is the solitary wave limit. The second, which arises
as the discriminant of a cubic which governs the normal form of the
linearized operator in a neighborhood of the origin, can also
be expressed in terms of the conserved quantities of the partial differential
equation and their derivatives with respect to the constants of integration
of the ordinary differential equation. This discriminant detects modulational
instabilities: bands of spectrum off of the imaginary axis which are connected
to the origin. As we have emphasized throughout this calculation can be considered to be a rigorous Whitham theory calculation.

This calculation hinges on the fact that the
underlying ordinary differential equation has
sufficient first integrals. As such it is doubtless
related to the multi-symplectic formalism of Bridges\cite{B}. As there are
a number of other equations for which the traveling wave ODE
has a integrable Hamiltonian formulation (Benjamin-Bona-Mahoney, NLS, etc)
one should be able to carry out the analogous calculation
in those cases. The additional structure provided by the scaling invariance
is also extremely helpful, as this allows one to simplify many of the
calculations but, at least in the Evan's function approach, it has not
been necessary.

We are somewhat puzzled by the fact that the Evan's function based
calculation gives a substantially simpler criteria for the existence
of a modulational instability than one based on a direct analysis of
the linearized operator. It must be true that the two discriminants
we've derived always have the same sign, as the predict the same phenomenon,
but we have been unable to
see this directly from the formulae. Often when apparently unconnected
quantities share a sign this sign has a topological or geometric
interpretation (for example as a Krein signature), so this may well
be the case here. Such an interpretation would be very interesting.

{\bf Acknowledgements:} The authors would like to thank Todd Kapitula for
several useful discussions, and Todd Kapitula and Mariana H{\v a}r{\v a}gu{\c s}
for sharing their numerical data. The insights provided were invaluable
in guiding the early stages of this work. J.C.B would like to acknowledge
support under NSF grants {\bf DMS-FRG}0354462  and {\bf DMS}0807584.

\section{Appendices}
\subsection{Identities}
In this section we derive a number of useful identities which will allow us
to relate various Jacobians which arise in the calculation. We define
the conserved quantities:
\begin{align*}
&T = \langle 1\rangle = 2\int_{u_-}^{u_+}\frac{du}{\sqrt{2 E + 2 a u - \frac{2}{p+2} u^{p+2} + c u^2}} \\
&M = \langle u\rangle = 2\int_{u_-}^{u_+}\frac{u du}{\sqrt{2 E + 2 a u - \frac{2}{p+2} u^{p+2} + c u^2}} \\
& P = \langle u^2\rangle = 2\int_{u_-}^{u_+}\frac{u^2 du}{\sqrt{2 E + 2 a u - \frac{2}{p+2} u^{p+2} + c u^2}} \\
& H =  \left\langle \frac{u_x^2}{2} - \frac{u^{p+2}}{p+2} \right\rangle =  2\int_{u_-}^{u_+}\frac{\left(E + a u - \frac{2 u^{p+2}}{p+2} + c u^2\right) du}{\sqrt{2 E + 2 a u - \frac{2}{p+2} u^{p+2} + c u^2}}.
\end{align*}
The classical action  $K = \oint u_x du=\int_0^Tu_x^2dx$ provides a useful generating
function, and is given by
\[
K = 2\int_{u_-}^{u_+}\sqrt{2 E + 2 a u - \frac{2}{p+2} u^{p+2} + c u^2}\;du.
\]
This integral has the advantage that it is regular at the endpoints and
can thus be differentiated in the form presented. It is obvious that the
derivatives are given by
\begin{align}
\frac{\partial K}{\partial a} &=  M\nonumber \\
\frac{\partial K}{\partial E} &=  T \nonumber\\
\frac{\partial K}{\partial c} &=  \frac{P}{2}.
\label{mixedpart1}
\end{align}
Note that these relations force certain relations among various $2\times 2$
Jacobians. For example we have
\[
\{M,P\}_{a,E} = \left\vert\begin{array}{cc}M_a & P_a \\ M_E & P_E\end{array}\right\vert =\left\vert\begin{array}{cc}K_{aa} & 2K_{ac} \\ K_{aE} & 2K_{cE}\end{array}\right\vert = 2\left\vert\begin{array}{cc}K_{aa} & K_{aE} \\ K_{ac} & K_{cE}\end{array}\right\vert = -2\{T,M\}_{a,c}.
\]
Similarly we have the relations
\begin{align*}
&\{T,M\}_{E,c} = -\frac{1}{2} \{T,P\}_{a,E} \\
&\{T,P\}_{a,c} =  \{M,P\}_{E,c}.
\end{align*}
There is another identity relating the gradients of $T$ and the conserved
quantities $M,P,H$ which is useful. Begin by noting that
\begin{align*}
\frac{1}{2} K +  \left< \frac{u^{p+2}}{p+2}\right> &= E T + a M + \frac{c}{2} P \\
\frac{1}{2} K - \left< \frac{u^{p+2}}{p+2}\right> &= H,
\end{align*}
where the former quantity is the Lagrangian of the traveling wave differential
equation. Adding these together, taking partial derivatives with respect to
$(a,E,c)$, and using the relations \eqref{mixedpart1} shows that
\[
\left(\begin{array}{cccc} T_a & M_a & P_a & H_a \\
  T_E & M_E & P_E & H_E \\
 T_c & M_c & P_c & H_c \end{array}\right) \left(\begin{array}{c} E \\ a \\ \frac{c}{2} \\ 1\end{array}\right) =  \left(\begin{array}{c}0 \\ 0 \\ 0 \\ 0 \end{array}\right)
\]
so that the gradients of any of $T,M,P,H$ can be expressed in terms of an
explicit linear combination of the other three. As noted in the text
the theory is most developed in terms of the first three quantities,
as these arise most naturally, but is stated in terms of the last three,
as these have the most natural physical interpretation.

There is another set of Jacobian identities which are useful. Differentiating
\eqref{dphi1} with respect to $E$ and subtracting $u_E$ times \eqref{dphi2}
gives the identity
\[
u_x u_{xE} - u_{xx} u_E = \{u,u_x\}_{x,E} = 1.
\]
Similarly we have the identities
\begin{align*}
& \{u,u_x\}_{x,a}=u \\
& \{u,u_x\}_{x,c}=\frac{1}{2}u^2.
\end{align*}
There are a number of other identities of this sort which can be derived in an analogous fashion.

\subsection{Analysis of $N(\mathcal{L}[u])$}

In this appendix, we give a detailed analysis of the null space $N(\mathcal{L}[u])$.  As above, we assume
$u$ is a solution of \eqref{gkdv} of period $T$ and define $\mathcal{H}=L^2(\RM / T\ZM;\CM)$.  Notice that from
Proposition \ref{Kernel} we know there that the functions $u_x$, $u_E$ both satisfy the differential equation
$\mathcal{L}[u]v=0$ when boundary conditions are ignored.  However, $u_E$ is not in general $T$-periodic due to
the variation in the period with respect to $E$.  This gives an indication that $u_x$ is the only
$T$-periodic solution of $\mathcal{L}[u]v=0$ unless one has $T_E=0$.  This observation
leads us to the main result of this appendix.

\begin{thm}
Considered as an operator on $\mathcal{H}$, $N(\mathcal{L}[u])=span\{u_x\}$.  In particular, up to constant multiples,
the equation $\mathcal{L}[u]v=0$ has only one solution in $\mathcal{H}$.
\end{thm}

\begin{proof}
Define $y_1(x)=\left(\frac{d u_-}{dE}\right)^{-1}u_E(x)$ and $y_2(x)=-V'(u_-)^{-1}u_x$ and note that
\begin{eqnarray*}
y_1(0)&=&1\;\;\;\;y_2(0)=0\\
y_1'(0)&=&0\;\;\;y_2'(0)=1.
\end{eqnarray*}
Moreover, $\mathcal{L}[u]y_j=0$ for $j=1,2$.  In the $y_1$, $y_2$ basis, an easy calculation then proves the monodromy
is expressed as
\[
m(0)=\left(
       \begin{array}{cc}
         1 & V'(u_-)\left(\frac{du_-}{dE}\right)^{-1}T_E \\
         0 & 1 \\
       \end{array}
     \right)
\]
and hence $0$ is a band edge of $\spec(\mathcal{L}[u])$.  It follows that there exists a second periodic element of $N(\mathcal{L}[u])$
if and only if $T_E=0$.  The proof that $T_E>0$ is a result of the following theorem by Schaaf, which completes the proof.
\end{proof}

\begin{lem}
Assume $G$ is a $C^3$ function on $(0,\infty)$ and that $G$ vanishes only at one point $x_0$ with
$G'(x_0)>0$.  Define
\[
A=\{x\in\RM : x<x_0\;\textrm{and } G(s)<0 \textrm{ for all
}s\in(x,x_0)\}
\]
and suppose for each $\alpha\in A$ there exists a periodic solution $x(t)>0$ of the
equation
\[
x''(t)+G(x(t))=0
\]
with initial data $x(0)=\alpha$, $x'(0)=0$.  Let $P(\alpha)$ denote the period of this solution
If $G$ satisfies the two conditions
\begin{eqnarray*}
G'(x)>0,\;x\in A&\Rightarrow& 5G''(x)^2-3G'(x)G'''(x)>0;\\
G'(x)=0,\;x\in A&\Rightarrow& G(x)G''(x)<0
\end{eqnarray*}
then $P$ is differentiable on $A$ and $P'(\alpha)>0$.
\end{lem}

In order to apply the above result in our case, define $G(x)=V'(x;a)$ and assume $G(x_a)=0$,
$G'(x_a)>0$, and $c=1$.  Define $A$ as above and notice that $G\in C^3(0,\infty)$.  For all $x$
such that $G'(x)>0$, we have $x^p>\frac{1}{p+1}$ and hence
\begin{eqnarray*}
(5G''^2-3G'G''')(x)&=&p(p+1)x^{p-2}\left[(p+1)(2p+3)x^p+3(p-1)\right]\\
&>&5p^2(p+1)x^{p-2}>0
\end{eqnarray*}
for such $x$. Also, if $G'(x)=0$, then $x^p=\frac{1}{p+1}$ and
hence at such points,
\begin{eqnarray*}
G(x)G''(x)=p\left[\frac{1}{p+1}-1-a(p+1)x^{p-1}\right]<0
\end{eqnarray*}
given that
\[
ax^{p-1}=a\left(\frac{1}{p+1}\right)^{\frac{p-1}{p}}>\frac{1}{(p+1)^2}-\frac{1}{p+1}
\]
for all $x\in A$. Hence, it follows for such $a$ that $T_E>0$ for
all periodic waves bounded by a homoclinic orbit in phase space.

Next, as pointed out in the text, the fact that the monodromy matrix at the origin is the identity plus
a rank one perturbation implies there is a linear combination of $u_a$ and
$u_E$ which is periodic.  This combined with the above lemma immediately implies
the following.

\begin{corr}\label{NullSpaces}
Considered as operators on $\mathcal{H}$, we have $N(\partial_x\mathcal{L}[u])=span\{\phi_0,\phi_1\}$ and
$N(\mathcal{L}[u]\partial_x)=span\{1,u\}$.
\end{corr}

We end our discussion with the following interesting remark.  By our above work we see that the requirement
$T_E\neq 0$ is sufficient of the origin to not be a double point for the sepctrum
of the Hill-operator $\mathcal{L}[u]$.  There is a geometric quantity which detects this same
information known as the \textit{Krein signature}.  For the operator $\mathcal{L}[u]$, the Krein
signature at the origin is easily shown to be $\tr(m'(0))$, where $m(0)$ is defined as above.
To see this, notice the characteristic polynomial for $m(\mu)$ can be expressed as
\[
\det[m(\mu)-\lambda I]=\lambda^2-\tr(m(\mu))\lambda+1
\]
which roots
\[
\lambda_{\pm}=\frac{\tr(m(\mu))\pm\sqrt{\tr(m(\mu))^2-4}}{2}.
\]
Thus, the solutions to the equation $\tr(m(\mu))=\pm2$ correspond to the periodic eigenvalues of the operator
$\mathcal{L}[u]$ and, moreover, $\mu$ is a double point of the periodic spectrum if and only if $\tr(m'(\mu))=0$.
From this discussion, it follows that $\tr(m'(0))$ must be a non-zero multiple of $T_E$, a fact which is proven in the following
proposition.

\begin{prop}
For the operator $\mathcal{L}[u]$, one has $sgn\left(\tr(m'(0)\right)=sgn(T_E)$.  As a result, $0$ is never a double point
of $\spec(\mathcal{L}[u])$ under the assumption $T_E\neq 0$.
\end{prop}

\begin{proof}
This proof is essentially given in Magnus and Winkler\cite{MW}.  Using
variation of parameters, one can
express $\frac{d}{d\mu}y_j$ and $\frac{d}{d\mu}y_j'$ in terms of the $y_j$ and $y_j'$.  Using the
facts that $\det(\MM(0))=1$ and $\tr(m(0))=2$, a bit of algebra eventually yields the expression
\[
\tr(m'(0))=sgn(y_1'(T))\int_0^T\left(\sqrt{|y_1'(T)|}y_2+sgn(y_1'(t))\frac{y_1(T)-y_2'(T)}{2\sqrt{|y_1'(T)|}}y_1\right)^2 dx
\]
It follows that
\begin{eqnarray*}
sgn(T(m'(0)))&=&sgn(y_1'(T))\\
&=&sgn\left(V'(u_-)\left(\frac{du_-}{dE}\right)^{-1}T_E\right)\\
&=&sgn(T_E)
\end{eqnarray*}
as claimed.
\end{proof}

\subsection{Negativity of $M_a$ }
In this appendix we show that $M_a(a,0)<0$ for $a$ sufficiently small and $c=1$.
Note that $M(a,0)$ can (after some rescaling) be expressed in the
form
\[
M(a,0)=\int_0^{r(a)} \frac{\sqrt{u}}{\sqrt{a + u - u^{p+2}}} du
\]
where $r(a)$ is the smallest positive root of $a + u - u^{p+2}=0$. This is
a smooth function of $a$ for $a$ small enough and satisfies
\[
r(a) = 1 + \frac{a}{p+1} + O(a^2)
\]
The main idea is to rescale the above so that the integral is over a
fixed domain and show that the integrand is a decreasing function of
$a$ on the new domain. Rescaling gives
\[
M(a,0)=\int_0^{1} \frac{\sqrt{u}}{\sqrt{\frac{a}{r^3} + \frac{u}{r^2} - r^{p-1}u^{p+2}}} du
\]
The quantity $\frac{a}{r^3} + \frac{u}{r^2} - r^{p-1}u^{p+2}$ satisfies
\[
\frac{a}{r^3} + \frac{u}{r^2} - r^{p-1}u^{p+2} = u - u^{p+2} + a(1 -\frac{2}{p+1}v - \frac{p-1}{p+1} v^{p+2})
\]
The second term is clearly positive on the open interval $(0,1)$, and thus
the rescaled integrand is a decreasing function of $a$, and thus $M_a<0$
for $a$ small enough.

\subsection{Evaluation of Virial-Type Identities}
In this section, we evaluate the two virial type inner products arising in the perturbation
analysis of section 4.3.  One of these is calculable for an arbitrary nonlinearity, while for the other
we must restrict to the case of power-law nonlinearities.  We proceed with the more general one first.

\begin{lem}\label{viral1}
$\left<\psi_0,L_1L_0^{-1}L_1\phi_1\right>=-T\{T,K\}_{a,E}$.
\end{lem}

\begin{proof}
Define an operator $\xi:\{g\in L^2(\RM/T\ZM):\left<g\right>\neq 0\}\to L^2(\RM/T\ZM)$ by
\[
\xi(g)=x-\frac{T}{\left<g\right>}\int_0^xg(s)ds.
\]
Then a straight forward computation shows that $L_0^{\dag}\xi(\phi_0)=f'(u)-c+\frac{T_E T}{\{T,\left<u\right>\}_{a,E}}$.
It follows that
\begin{eqnarray*}
\left<\psi_0,L_1L_0^{-1}L_1\phi_1\right>&=&2\{T,M\}_{a,E}\left<\left(f'(u)-c\right),L_0^{-1}u_{xxx}\right>\\
&=&2\{T,M\}_{a,E}\left<L_0^{\dag}\xi(\phi_0),L_0^{-1}u_{xxx}\right>\\
&=&T\left<\phi_0,u_{xx}\right>\\
&=&-T\{T,K\}_{a,E}
\end{eqnarray*}
as claimed.
\end{proof}

While the above expression holds for an arbitrary nonlinearity, we have found a closed form
expression of $\left<\psi_2,L_1L_0^{-1}L_0\phi_1\right>$ in the case of power non-linearities.
From the evaluation of the modulational instability index via Evans function techniques, it should
be that this inner product is calculable in the general case as well, although we have yet to be
able to do this.

\begin{lem}\label{viral2}
In the case of a power nonlinearity $f(x)=x^{p+1}$, we have
\begin{eqnarray*}
\left<\psi_2,L_1L_0^{-1}L_1\phi_1\right>&=&-T\{T,M\}_{E,c}\{T,K\}_{a,E}\\
&+&\frac{2-p}{p}\{T,M\}_{a,E}\left(M\{T,K\}_{a,E}-2\{T,M\}_{a,E}K\right)\\
&+&2c\{T,M\}_{a,E}\{T,M,K\}_{a,E,c}.
\end{eqnarray*}
\end{lem}

\begin{proof}
Notice that in the case of power-law nonlinearity, one has
\begin{eqnarray*}
\left<\psi_2,L_1L_0^{-1}L_1\phi_1\right>&=&-T\{T,M\}_{E,c}\{T,K\}_{a,E}\\
&+&\{T,M\}_{a,E}\left((2-p)\left<u^{p+1},L_0^{-1}L_1\phi_1\right>-2c\left<u,L_0^{-1}L_1\phi_1\right>\right),
\end{eqnarray*}
and hence we must evaluate $\left<u,L_0^{-1}L_1\phi_1\right>$ and $\left<u^{p+1},L_0^{-1}L_1\phi_1\right>$.
First, from the definition of $v_1$ in equation \eqref{v1} it follows that
\begin{eqnarray*}
\left<\psi_2,v_1\right>&=&\lambda_1\left<\psi_2,L_0^{-1}\phi_1\right>-\left<\psi_2,L_0^{-1}L_1\phi_1\right>\\
&=&-\frac{1}{2}\lambda_1\{T,M\}_{a,E}\{T,M,P\}_{a,E,c}-\{T,M\}_{a,E}\left<u,L_0^{-1}L_1\phi_1\right>
\end{eqnarray*}
Moreover, using the fact that $\psi_2=L_0^{\dag}\psi_1$ gives
\begin{eqnarray*}
\left<\psi_2,v_1\right>&=&\left<\psi_1,(\lambda_1-L_1)\phi_1\right>\\
&=&\lambda_1\left<\psi_1,\phi_1\right>+2\{T,M\}_{a,E}\left<\psi_1,u_{xxx}\right>\\
&=&-\frac{1}{2}\lambda_1\{T,M\}_{a,E}\{T,M,P\}_{a,E,c}-2\{T,M\}_{a,E}\left<\phi_2,u_{xx}\right>\\
&=&-\frac{1}{2}\lambda_1\{T,M\}_{a,E}\{T,M,P\}_{a,E,c}+\{T,M\}_{a,E}\{T,M,K\}_{a,E,c}
\end{eqnarray*}
and hence $\left<u,L_0^{-1}L_1\phi_1\right>=-\{T,M,K\}_{a,E,c}$.

Next, let the functional $\xi$ be as in Lemma \eqref{viral1} and notice that
\[
L_0^{\dag}\xi(u)=f'(u)-c-\frac{T}{M}\left(pu^{p+1}+a\right).
\]
It follows that
\begin{eqnarray*}
-\frac{Tp}{M}\left<u^{p+1},L_0^{-1}L_1\phi_1\right>&=&\left<L_0^{\dag}\xi(u)-(f'(u)-c),L_0^{-1}L_1\phi_1\right>\\
&=&\left<\xi(u),L_1\phi_1\right>+2\{T,M\}_{a,E}\left<(f'(u)-c),L_0^{-1}u_{xxx}\right>\\
&=&\frac{2T\{T,M\}_{a,E}T}{M}K-T\{T,K\}_{a,E}
\end{eqnarray*}
which completes the proof.
\end{proof}

\bibliography{GKdVbib}

\end{document}